\documentstyle[12pt,fullpage,amssym,saks,psfig,psfrag,epsfrag,here,
graphicx]{article}
\font\fFt=eusm10 scaled 1000

\newcommand{\K}{\,\mbox{\fFt K}}

\def\cc{\setcounter{equation}{0}   % THIS CLEARS THE COUNTER
\setcounter{figure}{0}\setcounter{table}{0}}

\begin{document}
\title{{\bf GEOMETRIC PROPERTIES OF QUASICONFORMAL MAPS AND SPECIAL
FUNCTIONS}}

\author{{\bf Matti Vuorinen}\\
Department of Mathematics, P.O. Box 4 (Yliopistonk. 5), \\
 FIN--00014 University of Helsinki, Finland\\[8mm]
%{\tt email: vuorinen@csc.fi}, \, \, FAX 358-9-19123213\\[10mm]
Lectures in the Vth Finnish-Polish-Ukrainian Summer School\\
in Complex Analysis, Lublin, 15-21 August 1996\\[10mm]
}

%\date{Preprint, \today \hskip 4mm ({\tt gmrv321.tex})}

\maketitle

\centerline{Published in: Bull. Soc. Sci. Lett. \L\'od\'z S\'er. Rech. D\'eform.  24  (1997), 7--58. }

\abstract
Our goal is to provide a survey of some topics in quasi\-con\-for\-mal analysis of
current interest. We try to emphasize ideas and leave proofs and
technicalities aside. Several easily stated open problems are given.
Most of the results are joint work with several coauthors. In
particular, we adopt results from the book authored by
Anderson-Vamanamurthy-Vuorinen \cite{AVV6}.

Part 1. Quasiconformal maps and
spheres

Part 2. Conformal invariants and special functions

Part 3. Recent results on special functions


\endabstract
%\bigskip
%\centerline{\tt DRAFT (FILE: lub9.tex)}
%\bigskip
%\footnotesize{20 Nov. 1996}

%\cc
\section{Quasiconformal maps and spheres}
Some current trends in multi-dimensional quasi\-con\-for\-mal analysis are
reviewed in \cite{G6}, \cite{G8}, \cite{I2}, \cite{V6}, \cite{V7}, \cite{Vu4}.
\bigskip
\begin{nonsec} \label{1.1} Categories of homeomorphisms
\end{nonsec}
Below we shall discuss homeomorphisms of a domain of
$\Bbb R^n$ onto another domain in
$\Bbb R^n,\ n\ge 2$. Con\-for\-mal maps provide a well-known subclass of
general homeomorphisms. By Riemann's mapping theorem this class is very
flexible and rich for $n=2$ whereas Liouville's theorem shows that, for
$n\ge 3,$ con\-for\-mal maps are the same as M\"obius transformations, i.e.,
their class is very narrow. Thus the unit ball
$B^n=\{x\in \Bbb R^n:\ |x|<1\}$ can be mapped con\-for\-mally only onto a
half-space or a ball if the dimension is $n \ge 3.$
Quasicon\-for\-mal maps constitute a convenient
interpolating category of maps, much wider than con\-for\-mal maps, and less
general than locally H\"older-continuous homeomorphisms. We also note
that bilipschitz maps are a subclass of quasi\-con\-for\-mal maps. Deferring
the definition of a quasisymmetric map to \ref{1.40}, we note that bilipschitz
maps are a subclass of quasisymmetric maps, which in turn are a subclass
of quasi\-con\-for\-mal maps.
%\end{nonsec}
\bigskip

\begin{nonsec} \label{1.2}  Modulus of a curve family
\end{nonsec}
Now follows perhaps the
most technical part of this paper, the definition of the
modulus of a curve family. The nonspecialist
reader may be relieved to hear that this
notion will be used later only in the definition of quasi\-con\-for\-mal
mappings and that an alternative definition of quasi\-con\-for\-mal mappings
can be given in terms of the geometric notion of linear dilatation (see 
\ref{1.8}). 
Let $G$ be a domain in $\Bbb R^n$ and let $\Gamma$ be a curve family in
$G$. For $p>1$ the $p$-{\it modulus} $M_p(\Gamma)$ is defined by
\begin{equation}
\label{eq1.3}
M_p(\Gamma)=\inf_{\rho\in F(\Gamma)}\int_G\rho^p dm\, ,
\end{equation}
where
$F(\Gamma)=\{\rho:G\to R\cup \{\infty\},\ \rho\ge 0$ Borel:
$\int_{\gamma}\rho ds\ge 1$ for all locally rectifiable $\gamma \in
\Gamma\}.$
\bigskip
The most important case is $p=n$ and we set $M(\Gamma)=M_n(\Gamma)$---in this case
we just call $M(\Gamma)$ the modulus of $\Gamma.$ 
The {\it extremal length} of
$\Gamma$ is $M(\Gamma)^{1/(1-n)}$. The modulus is {\it a con\-for\-mal 
invariant}, i.e. $M(\Gamma)=M(h\Gamma)$ if $h$ is a con\-for\-mal map
and $h\Gamma=\{h\circ \gamma:\ \gamma\in \Gamma\}$. For the basic
properties
of the modulus we refer the reader to \cite{V1}, \cite{Car}, \cite{Oh},
 \cite{Vu1}.

\begin{nonsec} \label{1.4}  Modulus and relative size
\end{nonsec}
For a domain $G\subset \Bbb R^n$ and $E,F\subset G$ denote
$$\Delta(E,F;G)=\{ {\mbox{all\, curves\, joining}}\ E\ 
{\mbox{ and}}\ F\ \mbox{in}\ G\}.$$
We define the {\it relative size} of the pair $E,F$ by
$$r(E,F)=\min\{d(E), d(F) \} /d(E,F)\, ,$$
 where $d(E) = \sup \{ |x-y|: x,y \in E \}$ and
$$d(E,F) =  \inf \{ |x-y|: x \in E, \mbox{and} \, y \in F \}.$$
If $E$ and $F$ are disjoint continua then $M(\Delta(E,F;\Bbb R^n))$
and $r(E,F)$ are simultaneously small or large. In fact, there are
increasing homeomorphisms $h_j: [0,\infty)\to [0,\infty)$ with
$h_j(0)=0,\ j=1,2,$ such that

\begin{equation} \label{eq1.5}
h_1(r(E,F))\le M(\Delta(E,F;\Bbb R^n))\le h_2(r(E,F))
\end{equation}
(see \cite{V1}, \cite{Vu1}). 
The explicit expressions for $h_j$ in \cite[7.41-7.42]{Vu1} involve special
functions.

\begin{nonsec} \label{1.6}  Quasiconformal maps
\end{nonsec} Let $K \ge 1. $ A homeomorphism $f:G\to G'$ is termed 
$K$-{\it quasi\-con\-for\-mal} 
if for all curve families $\Gamma$ in $G$
\begin{equation} \label{eq1.7}
M(f\Gamma)/K\le M(\Gamma)\le KM(f\Gamma).
\end{equation}
The least constant $K$ in (\ref{eq1.7}) is called {\it the maximal
dilatation } of $f.$

Note that con\-for\-mal invariance is embedded in this definition: for $K=1$
equality holds throughout in (\ref{eq1.7}).
\item{} This definition resembles the bilipschitz condition, but it
should be noted that quasi\-con\-for\-mal maps can transform distances in
a highly nonlinear and totally unlipschitz manner.

There are numerous
equivalent ways of characterizing quasi\-con\-for\-mal maps \cite{Car}. 
It often
happens that a mapping $K_1$-quasi\-con\-for\-mal in the sense of one
definition is $K_2$-quasi\-con\-for\-mal in the sense of another definition,
where $K_2$ depends from $K_1$ in an explicit way and, what is most
important,  $K_2\to 1$ if $K_1\to 1$. We shall next consider in
\ref{1.8}
an equivalent definition based on the linear dilatation.
We shall see
that in the case of this definition,
finding such a constant $K_2$ explicitly has
required a time span
as long as the history of higher-dimensional quasicon\-for\-mal maps.

\begin{nonsec} \label{1.8}  Linear dilatation
\end{nonsec}
For a homeomorphism
$f:G\to G',\ x_0\in G,\ r\in(0,d(x_0,\partial G))$, let
$$H(x_0,f,r)=\sup\{{{|f(x)-f(x_0)|}\over {|f(y)-f(x_0)|}}:\
|x-x_0|=|y-y_0|=r\},$$
$$H(x_0,f)=\limsup_{r\to 0} H(x_0,f,r).$$
Then $H(x_0,f)$ is called the {\it linear dilatation} of $f$ at $x_0$.
\vskip 4mm

\centerline{\psfig{figure=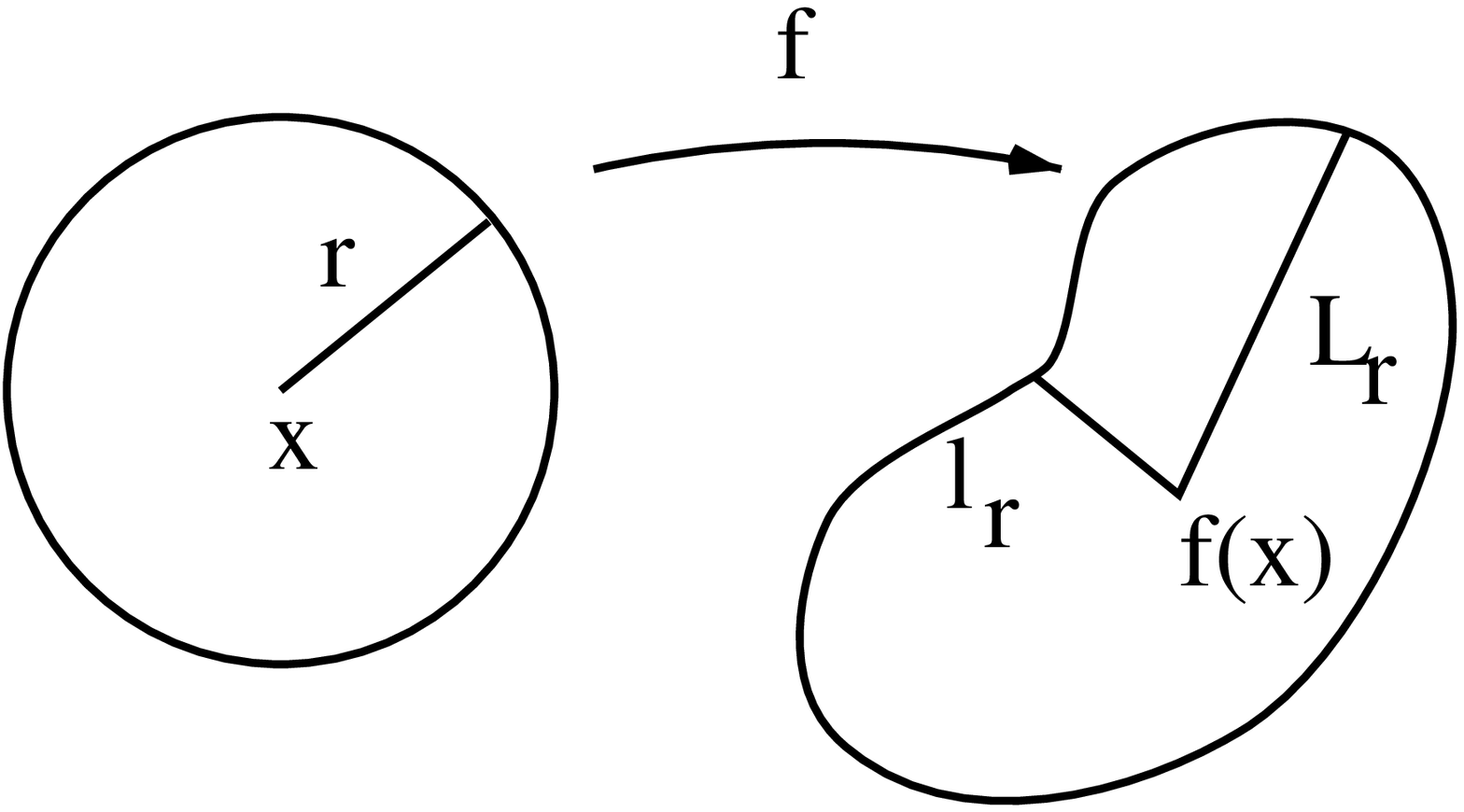,height=4cm}}

There is 
an alternative characterization of quasi\-con\-for\-mal maps,
 to the effect that a homeomorphism with  bounded linear dilatation
$$\sup \{ H(x,f) : x \in G \} \le L < \infty$$ is
quasi\-con\-for\-mal \cite{V1}. 
%In conclusion, the boundedness of linear dilatation is a
%characteristic property of quasi\-con\-for\-mal maps.
We shall next review the known estimates for the constant $L$ in
terms of the maximal dilatation.

Consider first the case $n=2 .$
A. Mori proved in \cite{Mor2} that if
 $f:G\to G', $ with $G, G' \subset \Bbb R^2,$ is $K$-quasi\-con\-for\-mal,
then for all $x_0 \in G$
\begin{equation} \label{Moribd}
H(x_0,f) \le e^{\pi K} .
\end{equation}
This bound is not sharp when $K \to 1 .$ The sharp bound
\begin{equation} \label{1.12}
H(x_0,f)\le \lambda(K)={u^2\over {1-u^2}}, \quad
 u=\varphi_K(1/\sqrt 2) \, ,
\end{equation}
is due to
Lehto, Virtanen, and
V\"ais\"al\"a  \cite{LVV} in the particular case $G ={\Bbb R}^2$ and due to
Shah Dao-Shing and Fan Le-Le \cite{SF} in the general case of 
a proper subdomain $G \subset {\Bbb R}^2 .$  
For the definition of the special function $ \varphi_K \equiv \varphi_{K,2},$
 see \ref{th2.26}.

Next we consider the case $ n \ge 2 .$
If $f:G\to G'\, , $ with $ G, G' \subset {\Bbb R}^n, $ 
is $K$-quasi\-con\-for\-mal then, by a 1962 result of
F.W. Gehring \cite[Lemma 8, pp. 371-372]{G1},
\begin{equation} \label{eq1.9}
H(x_0,f)\le d(n,K)\equiv \exp \left[ \left({K\omega_{n-1}\over 
\tau_n(1)}\right)^{1/(n-1)} \right]
\end{equation}
for all $x_0\in G,$
where $\omega_{n-1}=n\pi^{n/2}/\Gamma(1+{n\over 2})$ is the
$(n-1)$-dimensional surface area of the unit sphere $\partial B^n$,
 and $\tau_n$ is the
capacity of the Teichm\"uller condenser (see \ref{2.10}). For $n=2,$
the earlier result of A. Mori (\ref{Moribd}) is recovered as a particular case
of (\ref{eq1.9}), that is, $d(2,K) = e^{\pi K} .$
 Unfortunately
$d(n,K)\nrightarrow 1$ as $K\to 1$. In 1986 M. Vuorinen sharpened 
the bound (\ref{eq1.9}) to
\begin{equation} \label{1.10}
H(x_0,f) \le c(n,K) \equiv
 1+ \tau_n^{-1}(\tau_n(1)/K)<{1\over 10}d(n,K) \, .
\end{equation}
 Note that $c(n,K)\to 2$ as $K \to 1$ \cite[10.22, 10.32]{Vu1}.
In 1990 Vuorinen proved for a $K$-quasi\-con\-for\-mal map  
$f: \Bbb R^n \to \Bbb R^n$
of the whole space $\Bbb R^n$ \cite{Vu2} 
\begin{equation} \label{1.11}
H(0,f)\le \exp(6(K+1)^2\sqrt{K-1}) \equiv s(K)
\end{equation}
with the desirable property $s(K)\to 1$ as $K\to 1$.
In 1996 P. Seittenranta \cite{Se2} was able to prove a
similar
result for maps of proper subdomains $G$ of $\Bbb R^n:$
a $K$-quasi\-con\-for\-mal mapping $f:G \to G'$ satisfies
\begin{equation} \label{1.11b}
H(x_0,f) \le s(K)
\end{equation}
 for all
$x_0\in G$ with the same $s(K)$ as in (\ref{1.11}). Note that (\ref{1.11b})
would easily follow from (\ref{1.11}) if we could solve a local structure
problem stated below in \ref{2.50}. In fact,
slightly better bounds than (\ref{1.11}) and (\ref{1.11b}),
involving the special function $\tau_n$ are known.
Note also that for $n=2$ a sharper form of (\ref{1.11b}) holds by
(\ref{1.12})
and \cite{AVV2} since, for $K > 1,$

\begin{equation} \label{1.13}
 \exp(\pi(K-1)) \le \lambda(K)\le \exp(\pi(K-1/K)).
\end{equation}

\begin{nonsec} \label{1.18} Open problem
\end{nonsec} Can the upper bound (\ref{1.11b}) be replaced by $s(n,K)$
with $\lim_{n \to \infty} s(n,K) =1$ for each fixed $K >1$?

\begin{nonsec} \label{1.20}  Quasispheres and quasi\-circles
\end{nonsec}
If
$f:\Bbb R^n\to \Bbb R^n, \, n \ge 2,$ is $K$-quasi\-con\-for\-mal, then 
the set $fS^{n-1}$
is called a $K$-{\it quasisphere} or, if $n=2,$ a $K$-{\it quasicircle}.
 Here, as usual, $S^{n-1}=\partial B^n$ and
$B^n=\{x\in \Bbb R^n:\ |x|<1\}$.

Plane domains that are bounded by quasi\-circles, called quasidisks,
have been studied
extensively. See the surveys of Gehring \cite{G5}, \cite{G7}.
Compared to what is
known for the dimension $n=2$, very little is known in higher dimensions
$n\ge 3$. We shall formulate below some open problems, both for the
plane and the higher-dimensional case.

Part of the interest of quasispheres derives from the fact that these
sets can
have interesting geometric structure of fractal type. In fact, some of
the differences between the categories of bilipschitz and
quasi\-con\-for\-mal maps can be understood if one studies the geometric
structure of the images 
of spheres under these maps.

\begin{nonsec} \label{1.21}  Examples of quasi\-circles
\end{nonsec}
(1) Perhaps
the most widely known example of a nonrectifiable
quasi\-circle is the {\it snowflake curve}
(also called {\it von Koch curve}), which is constructed in the following way.
Take an equilateral triangle. To each side adjoin an equilateral
triangle whose base agrees with the middle-third segment of the side;
then remove this middle-third segment. Iterating this procedure
recursively ad infinitum we get as a result a nonrectifiable Jordan
curve of Hausdorff dimension $ > 1.$
Other similar examples are given in \cite{GV2}, \cite[p. 25]{G5}, 
and \cite[p. 110]{LV}.

%\hskip 2truecm
%\psfig{figure=lubpic3.ps,height=4cm}
\centerline{
\psfig{figure=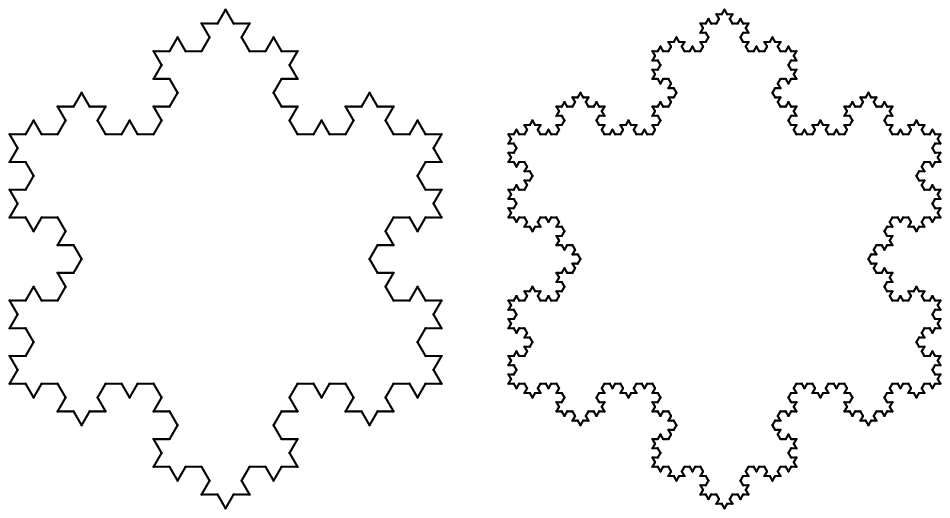,width=8.5cm}}
(2) The Julia set $J_f$ of an iteration $z \mapsto f(z)$ is the set of 
all those points that remain bounded under the repeated
iterations. As a rule, $J_f$ has an interesting fractal type structure,
and for suitable $f,$ $J_f$ is a quasicircle. For the case of
quadratic $f,$ see \cite{GM} and for rational $f$ see \cite{St}.

(3) Images of circles under bilipschitz maps are always rectifiable (and
hence of Hausdorff dimension 1) but they may fail to have tangents at
some points. In fact, bilipschitz maps are differentiable only almost
everywhere and if this \lq\lq bad set" of zero measure
is nonempty peculiar things may happen.
 See \cite{VVW} for a construction of a bilipschitz circle which is $(q,2)$-
thick in the sense of definition \ref{1.59} below.

(4) There are examples of Jordan domains with rectifiable boundaries which
are not bounded by quasispheres. For instance, the 
\lq\lq rooms and corridors"-type 
domains violating the Ahlfors condition in (\ref{eq1.31}) can be used.

(5) We next give a construction of a bilipschitz map 
$f : \Bbb R^2 \to \Bbb R^2$ with $f(0) =0$ which carries rays passing
through $0$ to \lq\lq logarithmic spirals" through $0.$ We first fix an
integer $p \ge 5$ and note that there exists $L \ge 1$ and an $L-$bilipschitz
mapping of the annulus $\overline{B}^2(p) \setminus B^2$ which is
identity on $S^1(p)$ and a restriction of the rotation 
$z \mapsto e^{i\theta} z, \theta \in (0,\pi/(2 p)),$ on $S^1(1).$
The boundary values of this map guarantee
that this mapping can be extended to an $L-$bilipschitz map of the
whole plane, which in the annuli $B(p^{k+1}) \setminus
\overline{B}(p^k)\, , k \in Z\, ,$
agrees with our original map up to conjugations by 
suitable rotations and dilations. For a similar construction, see
Luukkainen and V\"ais\"al\"a \cite[3.10 (4), 4.11]{LuV}.

(6) The univalent function
$$ f(z) = \int_0^z \mbox{exp} \{ ib \sum_{k=0}^{\infty} \zeta^{2^k} \} 
d \zeta, \quad b < {1 \over 4}\, ,$$
  defined in the unit disk $B^2,$
provides an analytic representation of a quasicircle $\Gamma = f(\partial D)$
that fails to have a tangent at each of its points. For details see
Ch. Pommerenke \cite[pp.304-305]{Po}.

\begin{nonsec} \label{1.22.} Particular classes of domains
\end{nonsec}
The unit ball in $\Bbb R^n$ is the standard domain for most applications
in quasi\-con\-for\-mal analysis. Since the early 1960's several classes of
domains have been introduced in studies on quasi\-con\-for\-mal maps.
It is not our goal to review such studies, but
we note that at least the following two types of domain 
classes have been studied:

(1) domains satisfying a geometric condition;

(2) domains characterized by conditions involving
   moduli of curve families, capacities, or other analytic conditions.

Domains of type (1) include so-called uniform domains and their various
generalizations. Domains of type (2) include, e.g., so-called
QED-domains. A domain $G \subset \Bbb R^n$ is called $c-$QED, $c \in (0,1]$
 if, for each pair of disjoint continua $F_1, F_2 \subset G,$
it is true that
$M(\Delta(F_1,F_2;G)) \ge c M(\Delta(F_1,F_2;\Bbb R^n)) .$
  There is a useful survey of some of these classes by
J. V\"ais\"al\"a \cite{V6}.

%It is not very clear to which class domains bounded by quasispheres
%should be included, perhaps they should be included in both classes
%(1) and (2).

Let us look at a property of the unit ball.
For nondegerate continua $E,F\subset B^n$ we have
$$M(\Delta(E,F;\Bbb R^n))\ge M(\Delta(E,F;B^n))\ge$$
$$M(\Delta(E,F;\Bbb R^n))/2\ge {1\over 2}h_1(r(E,F))\, $$
by \cite{G4} and (\ref{eq1.5}). (In particular, the unit ball is
$1/2$-QED.) For a domain $D \subset \Bbb R^n$ and $r_0>0$ we set
\begin{equation}
\label{eq1.23}
L(D,r_0)=\inf_{r(E,F)\ge r_0} M(\Delta(E,F;D)) \, ,
\end{equation}
where $E$ and $F$ are continua. For all dimensions $n\ge 2$
it is easy to construct \lq\lq rooms and corridors" type Jordan
domains with $L(D,r_0)=0$ (only simplest estimates of moduli are
needed from \cite[pp. 20-24]{V1}).
 For dimensions $n\ge 3$ one can construct such
domains also in the form 
$$D_g=\{(x,y,z)\in \Bbb R^3:\ x>0,\ |y|<g(x)\}$$
for a suitable homeomorphism $g:[0,\infty)\to [0,\infty),\ g(0)=0,\
g'(0)=0$; now the access to the \lq\lq ridge" 
$A \equiv \{(0,y,0):\ y\in
\Bbb R\}$ of the domain  gets narrower and 
narrower as we approach $A$ from within $D_g .$ 

\begin{figure}[H]
\begin{center}
\begin{psfrags}
    \psfrag{G1}[bl][bl]{\hspace{5mm} $\Gamma_1=\Delta (E,F;D)$}
    \includegraphics[width=2in]{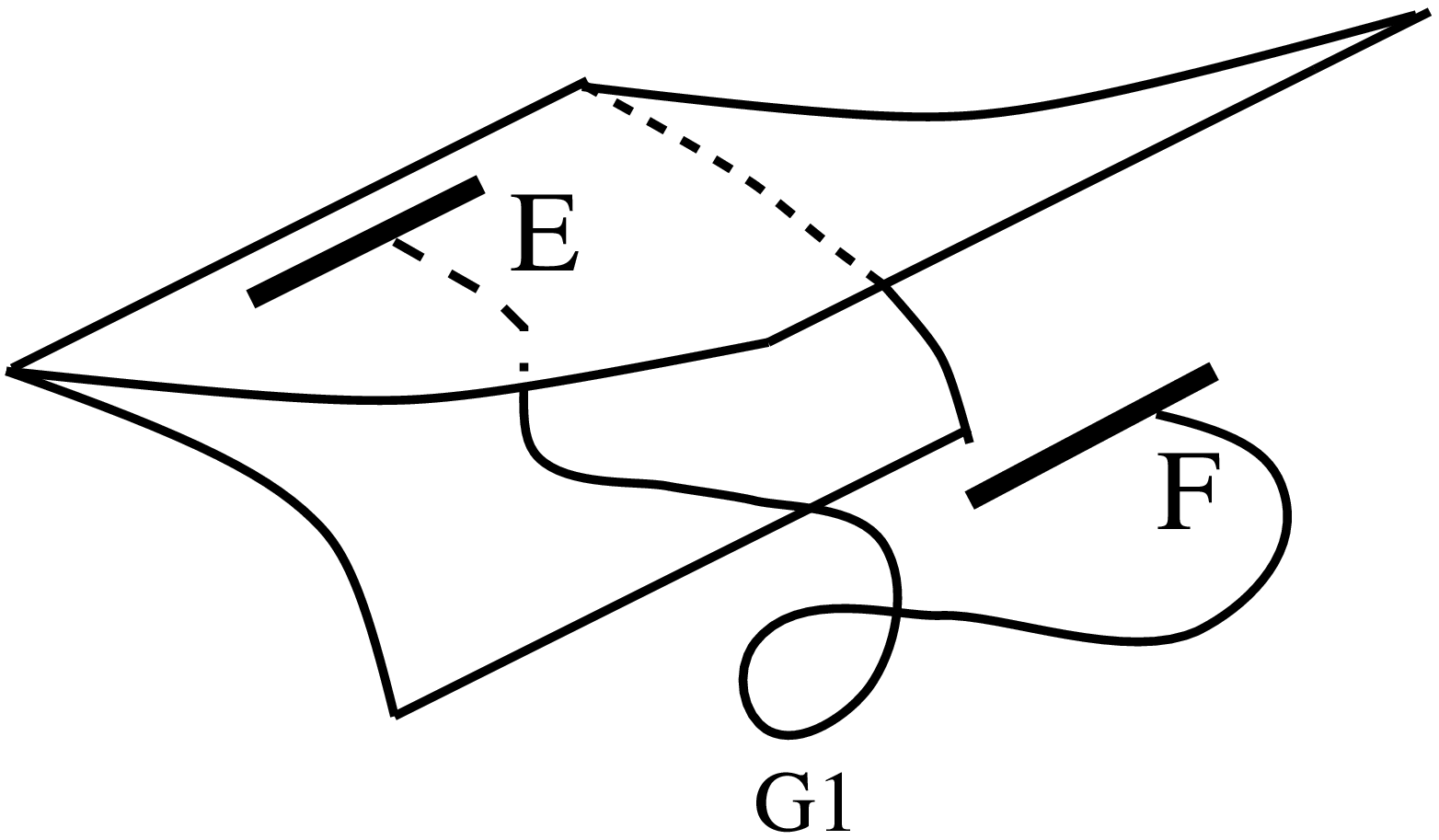}
\end{psfrags}
\end{center}
% \caption{Lublin picture 7a.}
\end{figure}

It is not difficult to show with the help of (\ref{eq1.5}) 
that the class of domains with
$L(D,r_0)>0$ is invariant under quasi\-con\-for\-mal maps of $\Bbb R^n$.
Hence we see that boundaries of domains with $L(D,r_0)=0$ cannot be
quasispheres.

\begin{figure}[H]
 \begin{center}
 \begin{psfrags}
     \psfrag{G2}[bl][bl]{$\Gamma_2$}
  \includegraphics[width=2in]{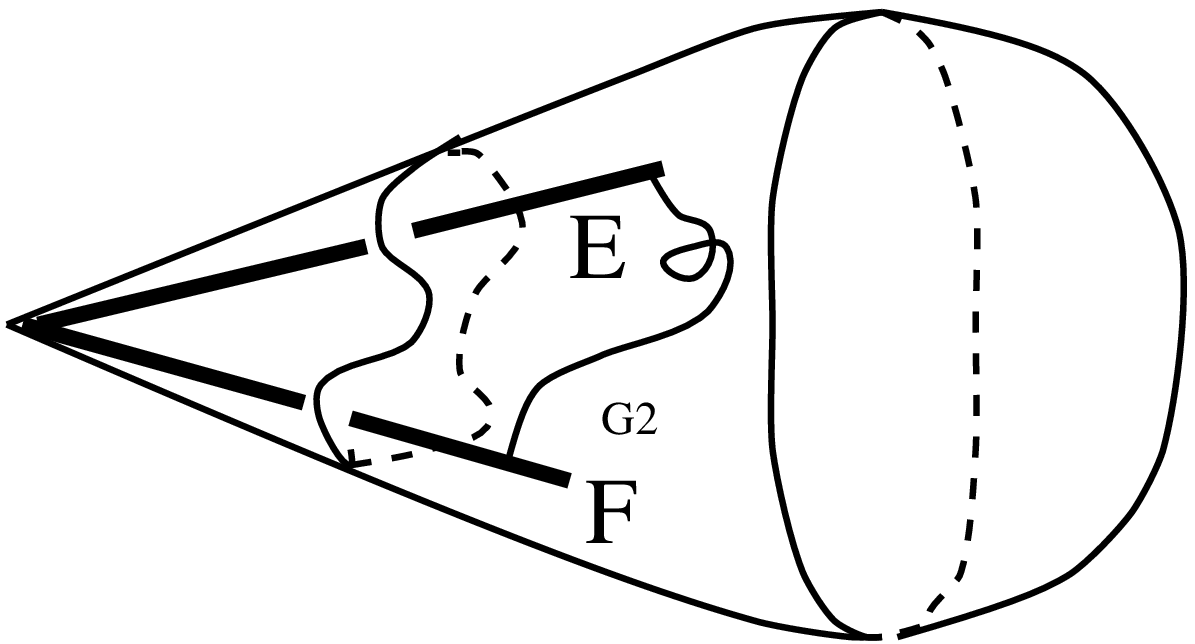}
 \end{psfrags}
 \end{center}
% \caption{Lublin picture 7b.}
 \end{figure}

%\vskip 0.5cm
%\centerline{\psfig{figure=lubp7a.eps,height=4cm}}
%\vskip 0.5cm

%\vskip 0.5cm
%\centerline{\psfig{figure=lubp7b.eps,height=4cm}}

One can also construct domains $D\subset \Bbb R^n$ such that for a pair
of disjoint continua $E,F\subset D$ with $r(E,F)=\infty$ we have
$M(\Delta(E,F;D))<\infty$.

%\vskip 0.5cm

\begin{nonsec} \label{1.24} Quasiconformal images of $B^3$
\end{nonsec}
By Liouville's theorem, the unit ball $B^n, \, n \ge 3,$ can be mapped con\-for\-mally
only onto another ball or a half-space. Gehring and V\"ais\"al\"a \cite{GV1}
created an extensive theory which gives necessary (and, in certain cases,
sufficient)
conditions for a domain to be of the form $fB^n$ where $f:B^n\to \Bbb R^n$
is quasi\-con\-for\-mal. They also exhibited several interesting domains
illuminating their results which we shall now discuss.

(1) The first example is an apple-shaped domain (cf. picture). By \cite{GV1}
such a domain cannot, in general, be mapped quasi\-con\-for\-mally onto $B^3$.

\centerline{\psfig{figure=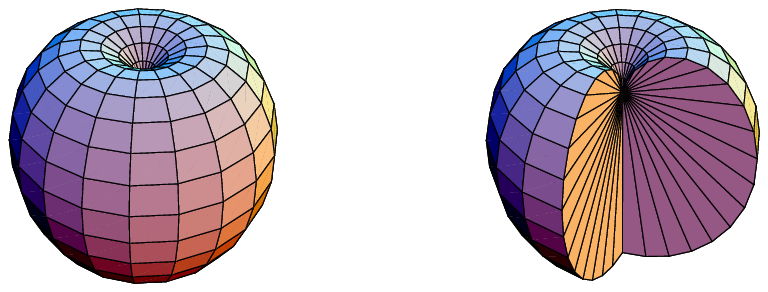,height=4cm}}

(2) On the other hand, there are onion-shaped domains that can be so mapped.

\centerline{\psfig{figure=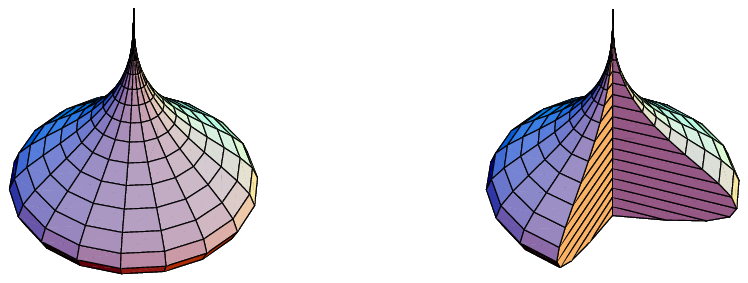,height=4cm}}

(3) In examples (1) and (2) above, the critical behavior takes place
near one boundary point at the tip of a spire.
In the case of an onion-shaped domain the spire is outwards-directed and for
apple-shaped domains it is inwards-directed. 
In this and the following example the
critical set consists of the edge of a boundary \lq\lq ridge". An example of a
domain with inward-directed ridge is shown (\lq\lq yoyo-domain") in the
picture below. The shape of the yoyo can be so chosen that the domain is
a quasi\-con\-for\-mal image of $B^3$.

\centerline{\psfig{figure=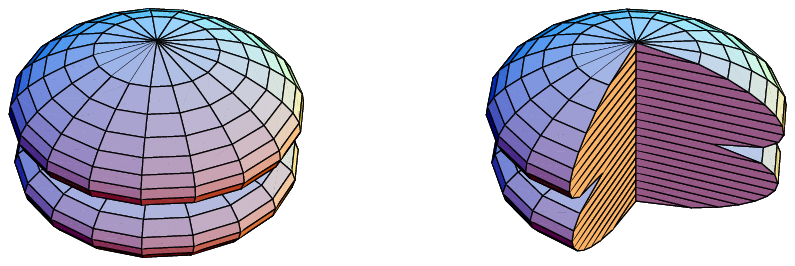,height=4cm}}

%\eject
(4) Consider now a \lq\lq ufo-shaped" domain where the ridge is
outward-directed (cf. the picture below).
 In this case the shape can be so chosen
that the domain is not quasi\-con\-for\-mally equivalent to $B^3$.

\centerline{\psfig{figure=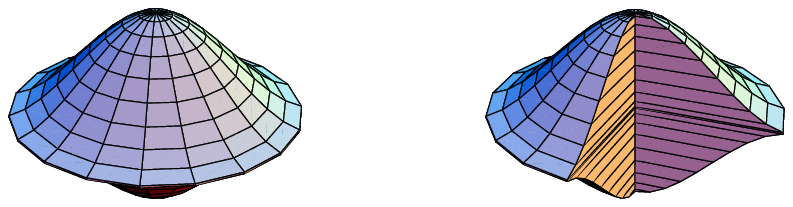,height=4cm}}

(5) P. Tukia \cite{Tu2} used an example of S. Rickman to construct a domain
whose boundary is the Cartesian product $K\times \Bbb R$  where $K$
is a snowflake-style curve with a periodic structure.
The domain underneath the surface fails to be quasi\-con\-for\-mally
equivalent to $B^3$. 

\centerline{\psfig{figure=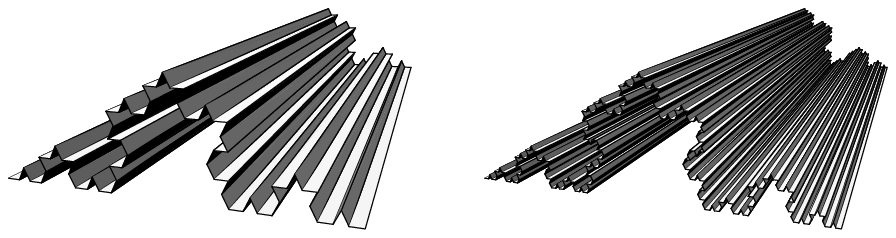,height=4cm}}

(6) Note that for dimensions $n\ge 3$ it is possible that a Jordan domain can
be quasi\-con\-for\-mally mapped onto $B^n$ but that
its complement fails to have this property.

\begin{nonsec} \label{1.30} Ahlfors' condition for quasi\-circles
\end{nonsec}
Quasicircles have been studied extensively and many characterizations
for them given by many authors. For interesting surveys, see \cite{G5}, \cite{G7}. 
Chronologically, one of the first characterizations
was given by L. V. Ahlfors in \cite{Ah1} and this result still continues to be
the most popular one and it reads as follows: A Jordan curve $C\subset
{\overline {\Bbb R}}^2$ is a quasi\-circle if and only if there exists a
constant $m \ge 1$ such that for all finite points $a,b\in C$
\begin{equation} \label{eq1.31}
\min\{d(C_1),d(C_2)\}\le m|a-b|\, ,
\end{equation}
where $C_1$ and $C_2$ are the components of $C\setminus \{a,b\}$ and
where $d$ stands for the Euclidean diameter.

Note that this formulation shows that (\ref{eq1.31}) 
guarantees the existence of
a $K$-quasi\-con\-for\-mal mapping 
$f:\Bbb R^2\to \Bbb R^2$ such that $C=fS^1.$
%We do not know the best upper bound for $K$ in terms of $m .$
 However, the least upper bound for $K$ in terms of $m,$
is not known.

%In any case we wish to call the reader's attention to the fact that
%showing that a Jordan curve is a $K$-quasi\-circle where $K$ is explicit,
%may require some effort. See also the open problems ...

\centerline{\psfig{figure=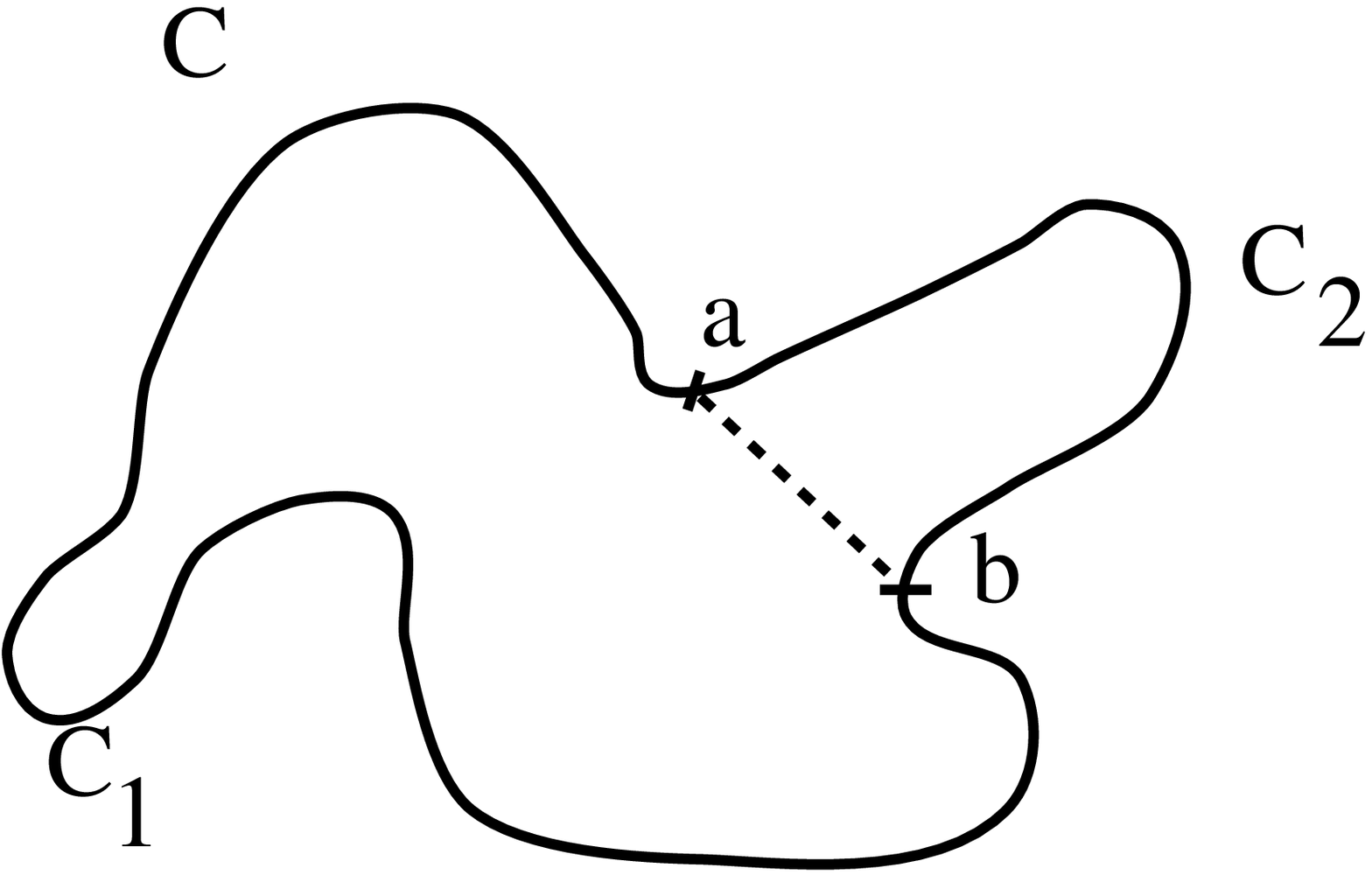,height=4cm}}

\begin{nonsec} \label{1.32} Open problem
\end{nonsec}
Generalize Ahlfors' condition to quasispheres.

\begin{nonsec} \label{1.33} Bilipschitz circles and spheres
\end{nonsec}
In harmony with our hierarchy of the categories of maps in 1.1, it is
natural to ask if a criterion similar to (\ref{eq1.31}) exists also for
bilipschitz circles or surfaces.
The case $n=2$ was settled by P. Tukia \cite{Tu1} in 1980 and also by
D. Jerison- C. Kenig \cite{JK} in 1982. 
The case $n\ge 3$ is open. Some results of this type were obtained
by S. Semmes \cite{S1}, \cite{S2} and T. Toro \cite{To1}, \cite{To2}.

\begin{nonsec} \label{1.34} Open problem
\end{nonsec}
Find the least $K$ for which a quadrilateral with given dimensions
is a $K$-quasi\-circle. A particular case is the rectangle.
R. K\"uhnau \cite[p. 104]{Kyuh2} has proved that a triangle with the
least angle $\alpha \pi (<  \pi/3)$
is a $K$-quasi\-circle with $K^2 \ge (1+d)/(1-d), d= |1-\alpha|,$ with
equality for the equilateral triangle $(\alpha =1/3 ) .$
(In fact, equality holds for all $\alpha \in (0,1/3)$ by S. Werner
\cite{We}.)
% Similar problems are unsolved also
%for quadrilaterals and other regions bounded by polygonal lines.

\begin{nonsec} \label{1.35} Open problem - triangle condition
\end{nonsec}
We say that a Jordan curve $C\subset \overline {\Bbb R}^2$ with
$\infty\in C$ satisfies a triangle condition if there exists a constant
$M\ge 1$ such that for all successive finite points $a,b,c\in C$ we have

\begin{equation} \label{eq1.36}
|a-b|+|b-c|\le M|a-c|
\end{equation}

Show that there exists a constant $K\ge 1$ such that $C=f\Bbb R$ where
$f:\Bbb R^2\to \Bbb R^2$ is $K$-quasi\-con\-for\-mal. Give $K=K(M)$ explicitly
in terms of $M$ with $K(M)\to 1$ as $M\to 1$.

\centerline{\psfig{figure=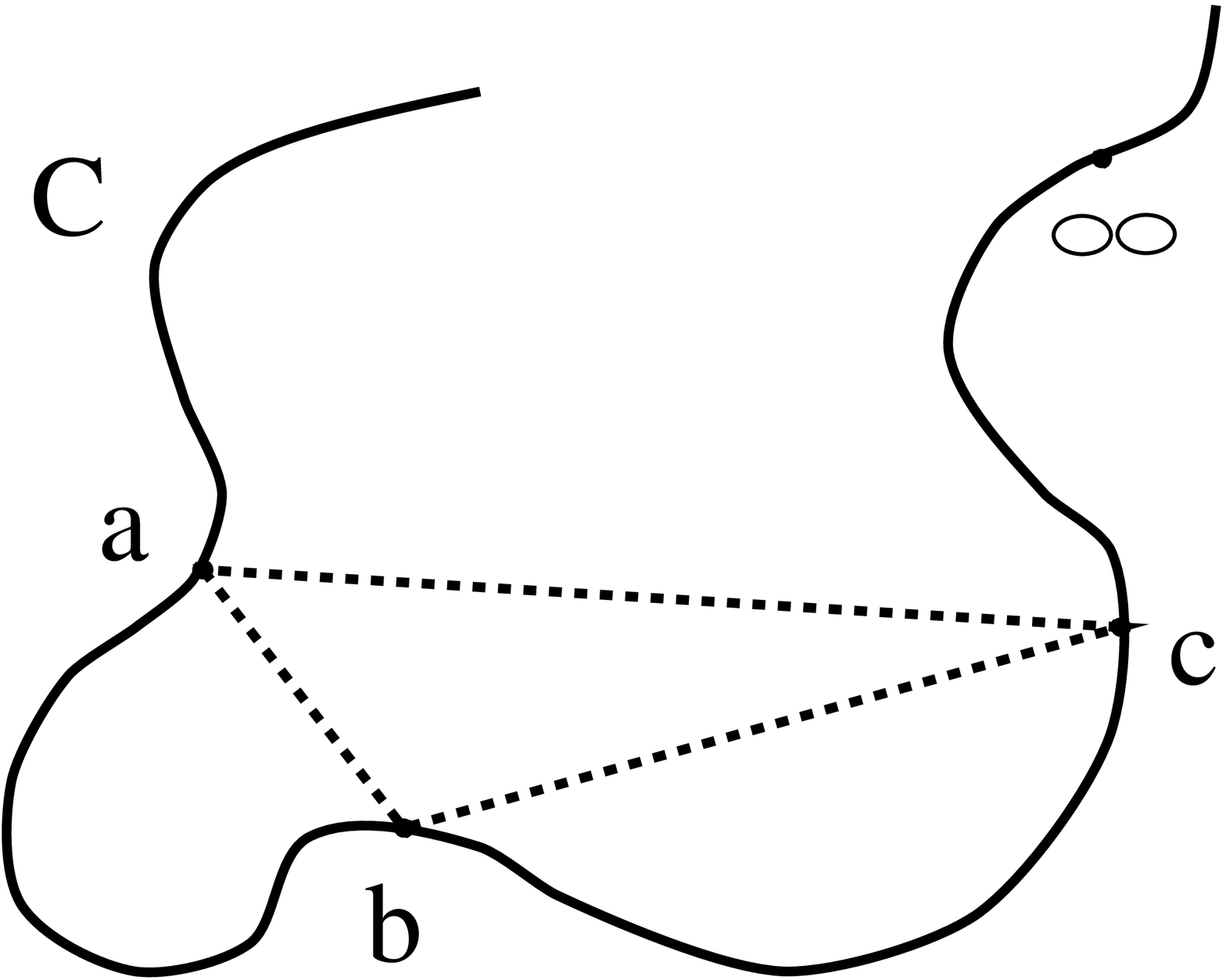,height=4cm}}

\begin{nonsec} \label{1.37} Remarks
\end{nonsec}
(1) From a result of S. Agard - F.W. Gehring \cite{AG} it follows that
$K(M)\ge 1+0.25(M-1)$ for $M\in (1,2)$.

(2) D. Trotsenko has informed the author (1996) about an idea to settle
the open problem \ref{1.35} with $K(M)\le 1+c_1\sqrt{M-1},\ c_1=10^5,$
for 
$M<1+10^{-13}$. See also \cite{Tr}.

\begin{nonsec} \label{1.40} Quasisymmetric maps
\end{nonsec}
Let $\eta:[0,\infty)\to [0,\infty)$ be a homeomorphism with $\eta(0)=0$
and let $f:G\to G'$ be a homeomorphism, where $G,G'\subset \Bbb R^n$.
We say \cite{TV1}  that $f$ is $\eta$-{\it quasisymmetric} if, for all
 $a,b,c\in G$ with
$a\not=c\, ,$
\begin{equation} \label {eq1.41}
{{|f(a)-f(b)|}\over {|f(a)-f(c)|}}\le \eta({|a-b|\over |a-c|})
\end{equation}

\begin{nonsec} \label{1.42}
Beurling - Ahlfors extension result
\end{nonsec}
A. Beurling and L.~Ahlfors \cite{BAh} introduced the class of
homeomorphisms $h : \Bbb R \to \Bbb R$ satisfying
\begin{equation} \label{eq1.55}
{1\over M}\le {{h(x+t)-h(x)}\over {h(x)-h(x-t)}}\le M
\end{equation}
for all $x \in \Bbb R, t > 0,$ and for some $M > 1.$  Such
homeomorphisms were later termed quasisymmetric. Note that, for maps of
the real axis, condition (\ref{eq1.55}) agrees with (\ref{eq1.41})
under the additional constraint $|a-b|=|a-c|.$
Beurling and Ahlfors also proved that a homeomorphism $f : \Bbb R \to \Bbb R$ of the
real axis can be extended to a $K$-quasi\-con\-for\-mal
map $f^*:\Bbb R^2\to \Bbb R^2$ iff $f$ satisfies (\ref{eq1.55}).
We remark that again there is a problem of finding the
optimal constant $K$ if $M >1$ is given. 
It is known by \cite[ p. 34 ]{L} that one
can choose $K \le \min \{ M^{3/2} , 2M -1 \} .$

\begin{nonsec} \label{1.43}
Quasisymmetry - quasi\-conformality
\end{nonsec}
If $f:G\to G'$ satisfies (\ref{eq1.41}) it follows easily that $H(x_0,f)\le
\eta(1)$ for all $x_0\in G$. By the alternative characterization of
quasi\-con\-for\-mality in terms of the linear dilatation \ref{1.8},
we thus see that quasisymmetric maps
constitute a subclass of quasi\-con\-for\-mal maps. As a rule, these two
classes of maps are different. However, if $G=\Bbb R^n$ then
quasi\-con\-for\-mal maps are $\eta$-quasisymmetric, by a result of P. Tukia
and J. V\"ais\"al\"a \cite{TV1}. 
Much more delicate is the question of finding for a given $K >1$
an explicit $\eta_K$
which is \lq\lq asymptotically sharp" when $K \to 1.$
 In \cite{Vu2} it was shown, for the first time, that
an explicit $\eta_{K,n}(t)$ exists which tends to $t$ as $K\to 1$:
If $f:\Bbb R^n\to \Bbb R^n, n \ge 2,$ is $K$-quasi\-con\-for\-mal, then
$f$ is $\eta_{K,n}$-quasisymmetric with

\begin{equation} \label{eq1.44}
\cases{\eta_{K,n}(1)\le \exp(6(K+1)^2\sqrt{K-1}),\cr
          \eta_{K,n}(t)\le \eta_{K,n}(1)\varphi_{K,n}(t),&$0<t<1$,\cr
          \eta_{K,n}(t)\le \eta_{K,n}(1)/\varphi_{1/K,n}(1/t),&$t>1.$\cr}
\end{equation}

Here $\varphi_{K,n}(t)$ is the distortion function in the quasi\-con\-for\-mal
Schwarz lemma (cf. Theorem \ref{th2.26}) with 

\begin{equation} \label{eq1.45}
\lambda_n^{1-\beta}r^\beta\le \varphi_{1/K,n}(r)\le \varphi_{K,n}(r)\le
\lambda_n^{1-\alpha}r^\alpha \, ,
\end{equation}
$\alpha=K^{1/(1-n)}=1/\beta,\  \lambda_n\in [4,2e^{n-1})$.
A $K$-quasi\-con\-for\-mal map of $B^n$ need not be quasisymmetric, but its
restriction to $\overline B^n(s),\ s\in (0,1)$, is quasisymmetric. In
fact, P. Seittenranta \cite{Se2} proved that for prescribed $K > 1 $
and $n\ge 2$, there exists an explicit $s\in (0,1)$ such that
$f|\overline B^n(s)$ is $\overline{\eta}_{K,n}$-quasisymmetric where 
 $\overline{\eta}_{K,n}$ is of the same
type as in (\ref{eq1.44}).

\begin{nonsec} \label{1.45} Linear approximation property
\end{nonsec}
Our examples of quasi\-circles in \ref{1.21} show that quasi\-circles need not
have tangents at any point. On the other hand, when $K\to 1$, we expect
that $K$-quasi\-circles become more like usual circles. We next introduce
a definition which enables us to quantify such a passage to the limit:

Given integers $n\ge 2,\ p\in\{1,...,n-1\}$, and positive numbers
$r_0>0,\ \delta\in(0,1)$, we say that a compact set $E\subset \Bbb R^n$
satisfies the {\it linear} {\it approximation} {\it property} with parameters
$(p,\delta,r_0)$ if for every  $x\in E$ and every $r\in (0,r_0)$ there
exists a $p$-dimensional hyperplane $V_r\ni x$ such that
$$E\cap B^n(x,r)\subset \{w\in \Bbb R^n:\ d(w,V_r)\le \delta r\}.$$

P. Mattila and M. Vuorinen proved in 1990 \cite{MatV}
that quasispheres satisfy this property.

\begin{theo} {} \label{1.45b}
 Let $K_2>1$ be such that
$$c=\eta_{K,n}(1)^{-2}/2>15/32$$
for all $K\in (1,K_2]$. Then a $K$-quasisphere $E=fS^{n-1}$ satisfies
the linear approximation property with parameters
\begin{equation}
\label{eq1.46}
(n-1,4g(K),d(E)g(K)),\ \mbox{where} \, g(K)=\sqrt{1-2c} \, .
\end{equation}
\end{theo}

Observe that here $\delta=4g(K)\to 0$ as $K\to 1$.

This limit behavior shows that, the closer
$K-1$ is to $0,$ the better $K$-quasispheres can be 
locally approximated by $(n-1)$-dimensional hyperplanes. Note 
that at a point $x\in E$ the approximating
hyperplanes $V_r$ may depend on $r$: they will very strongly depend on
$r$ if $x$ is a \lq\lq bad" point. An example of such bad behavior is a
quasi\-circle which logarithmically spirals in a neighborhood of a point
$x$.

%\centerline{\psfig{figure=lubpic13.ps,width=9cm}}

\centerline{\psfig{figure=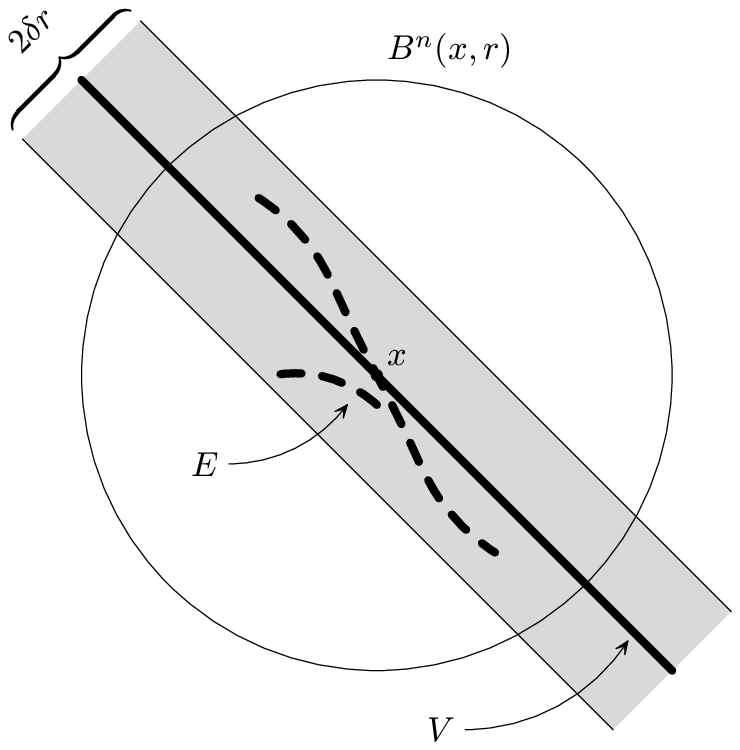,width=9cm}}

\begin{nonsec} \label{1.47} Jones' $\beta$-parameters
\end{nonsec}
In the same year as \cite{MatV} appeared, P. Jones \cite{Jo} introduced
\lq\lq $\beta$-parameters" for the analysis of geometric properties of plane
sets. In fact, the particular case $n=2,\ p=1,$ of the linear
approximation
property is very close to the
condition used by Jones in his investigations. Later on, Jones'
$\beta$-parameters were used extensively by  
C. Bishop - P. Jones \cite{BJ1}, G. David - S. Semmes
\cite{DS}, K. Okikiolu \cite{Ok}, and H. Pajot \cite{Paj}.

\begin{nonsec} \label{1.48} Open problem
\end{nonsec}
For $n=2$ the parameter $\delta$ of the linear approxi\-mation property in
(\ref{eq1.46}) is roughly $\sqrt{K-1}$. Can this be reduced, say to 
$K-1,$ when  $K$ is close to $1$?

\begin{nonsec} \label{1.49} Open problem
\end{nonsec}
The Hausdorff dimension of a $K$-quasicircle has a majorant of
the form $1 + 10(K-1)^2$ (see \cite{BP2}, \cite[1.8]{MatV}). 
Is there a similar bound for the Hausdorff dimension 
of a $K$-quasisphere in $\Bbb R^n,$ e.g. in the form
$n-1+ c(K-1)^2$ where $c$ is a constant?

\begin{nonsec} \label{1.50} Rectifiability of quasispheres
\end{nonsec}
Snowflake-type quasi\-circles provide examples of locally nonrectifiable
curves. We now briefly review conditions under which quasi\-circles will
be rectifiable. If $f:\Bbb R^n\to \Bbb R^n$ is $K$-quasi\-con\-for\-mal and
$t\in (0,1/2)$, then for convenience of notation we set

\begin{equation} \label{eq1.51}
K(t)=K(f|A(t)),\ A(t)=\bigcup_{x\in S^{n-1}}B^n(x,t)\, .
\end{equation}

\centerline{\psfig{figure=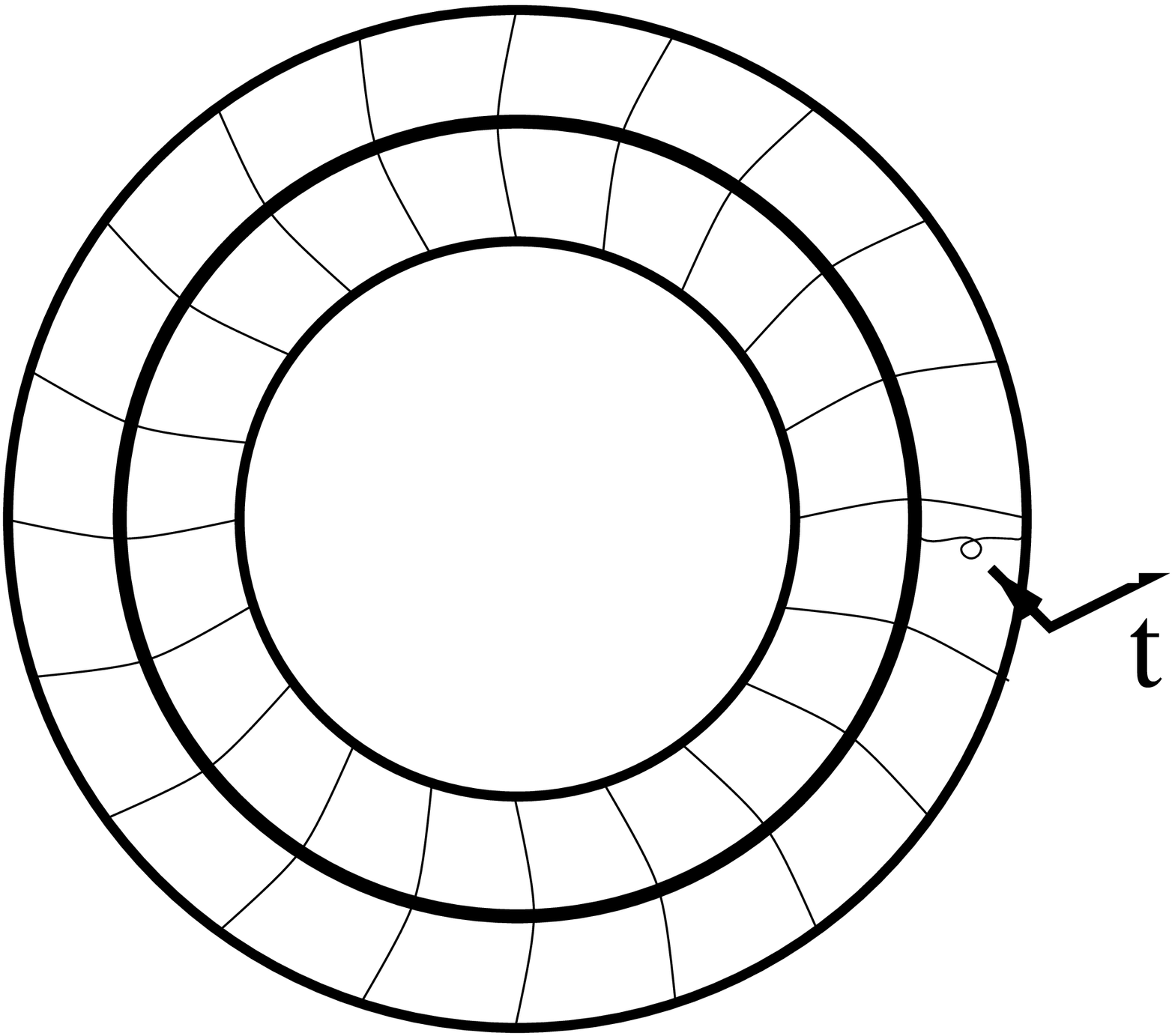,height=4cm}}
A natural question is this: Does $K(t)\to 1$ as $t \to 0 $
imply that $fS^{n-1}$
is rectifiable? For $n=2,$ J.\ Becker and Ch.\ Pommerenke \cite{BP1} have shown
that the answer is in the negative. Imposing a stronger condition for
the convergence $K(t)\to 1,$ we have a positive result \cite{MatV}:

\begin{theo}{} \label{thMatV}
If
\begin{equation} \label{eq1.52}
\int_0^{1/2}{{1-\alpha(t)}\over t}dt<\infty,\ \alpha(t)=K(t)^{1/(1-n)}\, ,
\end{equation}
then $fS^{n-1}$ is rectifiable.
\end{theo}

 An alternative proof of Theorem \ref{thMatV} was given by
 Yu. G. Reshetnyak in \cite[pp. 378-382]{Re2}. For some related results
see also \cite{GuV}. For
$n=2$ one can replace condition (\ref{eq1.52}) by a slightly weaker one,
as shown in \cite{ABL}, \cite{Carle}.

\begin{nonsec} \label{1.53} Quasiconformal maps of $S^{n-1}$
\end{nonsec}
Many of the peculiarities of quasi\-con\-for\-mal maps exhibited above are
connected with the interesting geometric structure of quasispheres.
We will now briefly discuss the simplest case when $f:\Bbb R^n\to
\Bbb R^n$ is a $K$-quasi\-con\-for\-mal map with $fS^{n-1}=S^{n-1}$. Let
$g=f|S^{n-1}$. Then $H(x,g)\le H(x,f)$ for every $x\in S^{n-1}$.
By the alternative characterization mentioned in \ref{1.8},
 we see that if $n-1\ge 2,$ then $g$ is
quasi\-con\-for\-mal
[note: we have not defined quasi\-conformality in dimension 1].
Thus for $n\ge 3$ the restriction $g$ satisfies all the properties of a
quasi\-con\-for\-mal map. In particular, $g$ is absolutely continuous with
respect to $(n-1)$-dimensional Hausdorff measure on $S^{n-1}$. For
$n=2$ the situation is drastically different, as the
following result of Beurling and Ahlfors shows.

\begin{nonsec} \label{1.54} Beurling - Ahlfors' singular function
\end{nonsec}
In \cite{BAh} Beurling and Ahlfors constructed a homeomorphism $h:\Bbb R\to
\Bbb R$ satisfying the condition (\ref{eq1.55})
%\begin{equation} \label{eq1.55}
%{1\over M}\le {{h(x+t)-h(x)}\over {h(x)-h(x-t)}}\le M
%\end{equation}
for some $M>1$ such that $h$ is not absolutely continuous with respect to
1-dimensional Lebesgue measure. By their extension result mentioned
in \ref{1.42}, $h$ is the restriction of a quasi\-con\-for\-mal mapping $h^*$ of
$\Bbb R^2$. If $g$ is a M\"obius transformation with $g(S^1)=\Bbb R$, then
the conjugation $g^{-1}\circ h^*\circ g$ is the required counterexample.

\begin{nonsec} \label{1.56} Tukia's quasisymmetric function
\end{nonsec}
Answering a question of W.K. Hayman and A. Hinkkanen, P. Tukia
constructed in \cite{Tu3} an example showing that 
a quasisymmetric  map $f$ of $\Bbb
R$ can map a set $E,$ with  H-dim $E<\varepsilon$ 
onto a set with $\mbox{H-dim}
({\Bbb R} \setminus fE) < \varepsilon$. See also \cite{BS} and \cite{Ro}.

\begin{nonsec} \label{1.57} Thick sets
\end{nonsec}
We conclude this section with a discussion of a property opposite to
the linear approximation property. Let $c>0,\ p\in \Bbb N$. We say that
$A\subset \Bbb R^n$ is $(c,p)$-{\it thick} if, for every $x\in A$ and for all
$r\in(0,d(A)/3),$ there exists a $p$-simplex $\Delta$ with vertices in
$A\cap B^n(x,r)$ with $m_p(\Delta)\ge cr^p$ \cite{VVW}, \cite{V5}.

Snowflake-type curves are examples of $(c,2)$-thick curves. One can even
show that for every $K>1$ there are $({\sqrt K-1\over 768},2)$-thick
$K$-quasi\-circles. For this purpose one uses a snowflake-style
construction, but replaces the angles $\pi\over 3$ by smaller ones that
tend to $0$ as $K\to 1$ \cite{VVW}.

A condition similar to thickness is the notion of wiggly sets \cite{BJ2}.

\begin{nonsec} \label{1.58}
Open problem
\end{nonsec}
Are there quasispheres in $\Bbb R^n,\ n\ge 3$, which are $(c,n)$-thick
for some $c>0$?

\begin{nonsec} \label{1.59}
Open problem
\end{nonsec}
Let $f:\Bbb R^n\to \Bbb R^n$ be a homeomorphism and let $fS^{n-1}$ be
$(c_1,n)$-thick. Is it true that
$$\mbox{H-dim}(fS^{n-1})\ge n-1+c_2>n-1?$$

For $n=2$ the answer is known to be in the affirmative \cite{BJ2}.

\begin{nonsec} \label{1.60} Additional references
\end{nonsec} The change of Hausdorff dimension under quasicon\-for\-mal maps
has been studied recently in \cite{IM2} and \cite{Ast}. A subclass of 
quasi\-circles, so-called {\it asymptotically con\-for\-mal curves}, has been
studied, for instance, in \cite{BP1}, \cite{ABL}, \cite{GR}.

\begin{nonsec} \label{1.61} More open problems
\end{nonsec} Some open problems can be found in 
\cite[p. 193]{Vu1}, \cite{AVV3}.

\begin{nonsec} \label{1.62} Books
\end{nonsec} The existing books on quasi\-con\-for\-mal maps include
\cite{Car}, \cite{KK}, \cite{L}, \cite{LV}, \cite{V1}. Generalizations
to the case of noninjective mappings, so-called {\it quasiregular mappings}, 
are studied in \cite{HKM}, \cite{I1}, \cite{IM2},
\cite{Re2}, \cite{Ri}, \cite{V2}, \cite{Vu1}.
%\eject
%SECTION 2 BEGIN
\cc

\section{Conformal invariants and special functions}

In this section we try to answer some fundamental questions such as:

a. Why are conformal invariants used in geometric function theory?

b. Why are special functions important for conformal invariants?

c. What are some of the open problems of the field?

In what follows we will provide some answers to these questions, as well
as pointers to the literature for further information. In a nutshell
our answer to \lq\lq a" is provided by the developments in geometry and
analysis that emerged from Klein's Erlangen Program and to \lq\lq b" by the
fact that the solution to some con\-for\-mally invariant extremal problems
involve special functions.

\begin{nonsec} \label{2.1} Klein's Erlangen Program
\end{nonsec}
The genesis of F. Klein's Erlangen Program is attached usually
to the year 1872 when Klein became a professor at the University of
Erlangen. In this program, the idea of using group theory to study
geometry was crystallized into a form where the following conceptions
played a crucial role

- use of isometries (\lq\lq rigid motions") and invariants

- two configurations are regarded equivalent if one can be carried to
the other by a rigid motion (group element)

- the basic \lq\lq models" of geometry are

(a) Euclidean geometry

(b) hyperbolic geometry (Bolyai-Lobachevskii)

(c) spherical geometry

The main examples of rigid motions are provided by various subgroups
of M\"obius transformations of 
$\ \overline{\Bbb R}^n=\Bbb R^n\cup\{\infty\} .$
 The group of M\"obius
transformations is generated by reflections in $(n-1)$-dimensional
spheres and hyperplanes.

\begin{nonsec} \label{2.2}
Geometric invariants
\end{nonsec}
In each of the models of Klein's geometries, there are natural metrics
that are invariant under \lq\lq rigid motions". For spherical geometry,
such a metric is the chordal metric, defined in terms of the
stereographic projection $\pi:\overline {\Bbb R}^n\to S^{n-1}({1\over
2}e_n,{1\over 2})$ by
$$q(a,b)={|a-b|\over {\sqrt{1+|a|^2}\sqrt{1+|b|^2}}}=|\pi a-\pi b| \, ,$$
$$  \pi x = e_{n+1}+(x-e_{n+1})/|x-e_{n+1}|^2 , $$
for $a,b, x \in \Bbb R^n$. The absolute (cross) ratio is defined as follows
$$|a,b,c,d|={q(a,c)q(b,d)\over q(a,b)q(c,d)}={|a-c||b-d|\over
|a-b||c-d|}.$$
Its most important property is invariance under M\"obius transformations.

We shall next consider a few examples of geometric invariants in the
sense of Klein.

\begin{nonsec} \label{2.3} Hyperbolic geometry
\end{nonsec}
For distinct points $a,b\in B^n$ 
let $a_*,\ b_*\in \partial B^n$ be distinct points such that the
quadruple $a_*,a,b,b_*$ can be moved by a rigid motion $T$ (=M\"obius
selfmap of $B^n$) to $(-e_1,0,\lambda e_1,e_1),\ \lambda \in (0,1).$ Then
$T^{-1}(-e_1,e_1)$ is an arc of an orthogonal circle through $a$ and $b$.
We define the {\it hyperbolic} {\it metric} $\rho$ by
\begin{equation} \label{eq2.4}
\rho(a,b)=\log|a_*,a,b,b_*|.
\end{equation}
By M\"obius invariance of the absolute ratio we see that $\rho$ is
invariant under M\"obius selfmaps of $B^n$.

In addition to a metric, another fundamental notion of hyperbolic
geometry is the hyperbolic volume of a polyhedron.
For $n=2$ and $a,b,c\in B^2$ let $\alpha,\beta,\gamma$ be the angles of
a triangle with vertices $a,b,c$. Then the {\it hyperbolic} {\it area}
of the triangle is
\begin{equation} \label{eq2.5}
v(a,b,c)=\pi-(\alpha+\beta+\gamma).
\end{equation}
Also hyperbolic area is invariant under M\"obius selfmaps of $B^2$.

A natural question is whether it is possible to define similar
geometries in domains not homeomorphic to $B^n$. In every proper
subdomain of $\Bbb R^n$ one can define the {\it absolute} {\it ratio}
{\it metric} by
\begin{equation} \label{eq2.6}
\cases{\delta_G(a,b)=\log(1+r_G(a,b))\cr
        r_G(a,b)=\sup\{|a,c,b,d|: c,d\in \partial G\}\cr}
\end{equation}
(see \cite{Se1}). Clearly this is a M\"obius-invariant metric, and it can be
shown that $\delta_{B^n}\equiv \rho$ for $G=B^n .$  Some of the basic
properties of $\delta_G$ are proved in \cite{Se1}. Another metric, 
the so-called {\it Apollonian} {\it metric},  defined by
\begin{equation} \label{eq2.7}
\alpha_G(a,b)=\sup\{\log|c,a,b,d|:c,d\in \partial G\}
\end{equation}
was studied recently in \cite{Be} and \cite{Se1} (strictly speaking, 
$\alpha_G$ is only a pseudometric).

\begin{nonsec} \label{2.8}
Origin of quasi\-conformal maps
\end{nonsec}
Klein's Erlangen Program received wide acclaim, and similar ideas proved
fruitful also in geometric function theory. H.A. Schwarz (Schwarz
lemma), H. Poincar\'e, and C. Carath\'eodory were some of the eminent
promoters of these ideas.

It is in this stage of the mathematical evolution that H. Gr\"otzsch wrote
his now famous 1928 paper which was to become the first paper on plane
quasi\-con\-for\-mal maps. It is
sometimes pointed out that, a century earlier in his theory of surfaces,
Gauss had studied notions that were close to quasiconformal maps.
%Sometimes it is pointed out that a century earlier
%Gauss studied in his theory of surfaces notions close to quasi\-con\-for\-mal
%maps. 
One of the important tools introduced by Gr\"otzsch was a new
con\-for\-mal invariant, the modulus of a quadrilateral.
Remarkable progress took place in 1950 when  L. Ahlfors and A. Beurling
found a new con\-for\-mal invariant, the extremal length of a curve family
(cf. 1.1) which soon became a popular tool in geometric function theory
(\cite{G6}, \cite{J1}, \cite{Kuz}, \cite{Rod}). 
Higher-dimensional quasi\-con\-for\-mal maps entered the stage
first in a note by M.A. Lavrentiev in 1938 but the systematic study was
started only in 1959 by C. Loewner, F.W. Gehring, B. Shabat, and
J. V\"ais\"al\"a.

\begin{nonsec} \label{2.8b}
Liouville's theorem
\end{nonsec}
Soon after the publication of Riemann's famous mapping theorem concerning
conformal maps of simply-connected
plane domains, Liouville proved that, in striking
contrast to the two-dimensional case, the only $C^3$ conformal maps of
subdomains of $\Bbb R^n, n \ge 3 ,$ are restrictions of M\"obius 
transformations. Under weaker differentiability hypotheses this result was
proved by F. W. Gehring \cite{G2} and Yu. G. Reshetnyak \cite{Re1}.
See also B. Bojarski and T. Iwaniec \cite{BI1}.
Yu. G. Reshetnyak has created so-called {\it stability theory,}
which is a study of properties of $K$-quasiconformal and
$K$-quasiregular maps with small $K-1.$ The main goal of this theory is
to find quantitative ways to measure the distance of these mapping classes
from M\"obius maps. The fundaments of this theory
are presented in \cite{Re2}. In spite of the many results in \cite{Re2},
some very basic questions are still open, see \ref{2.9b} below.
Significant results on stability theory were proved by
V. I. Semenov \cite{Sem1}, \cite{Sem2} and others.  

\begin{nonsec} \label{2.9}
Main problem of quasi\-conformal mapping theory
\end{nonsec}
We recall from Section 1 the definition of a $K$-quasi\-con\-for\-mal
map $f:G\to G'$ where $G$ and $G'$ are domains in $\overline {\Bbb
R}^n:$ a homeomorphism
$f$ is $K$-quasi\-con\-for\-mal if, for all curve families $\Gamma$ in $G',$
$$(*)\ \  M(f\Gamma)/K\le M(\Gamma)\le KM(f\Gamma).$$
A main problem of quasi\-con\-for\-mal mapping theory in $\Bbb R^n$
is how to extract explicit \lq\lq geometric information" from (*)
preserving \lq\lq asymptotic sharpness" as $K\to 1.$

It should be noted that the vast majority of results on quasi\-con\-for\-mal
mappings are not sharp in this sense. Of the results below, the Schwarz
lemma \ref{th2.26} is an example of an asymptotically sharp result.

We next list three simple ideas that might be used when studying this
problem.

{\bf Idea 1.} Use \lq\lq canonical situations", where the modulus of $\Gamma$
can be computed explicitly, as comparison functions.

{\bf Idea 2.} Idea 1 and the basic inequality (*) lead to (nonlinear)
constraints which we need to simplify.

{\bf Idea 3.} Try to relate (*) to \lq\lq geometric notions" distances,
metrics, etc. This leads to con\-for\-mally invariant extremal problems,
whose solutions can often be expressed in terms of {\it special}
{\it functions}.

\begin{nonsec} \label{2.9b} Open problem
\end{nonsec}
Let $f : {\Bbb R}^n \to {\Bbb R}^n $ be $K-$quasiconformal map normalized by
 $f(0) = 0, \, f(e_1) =e_1$ and let $I$ stand for the class of all
isometries
$h$ of ${\Bbb R}^n$ with $h(0) = 0, \, h(e_1) = e_1 . $ Find an explicit
and concrete upper bound for
$$ \varepsilon(K,n) \equiv 
\inf_{A \in I} \sup \{ |f(x) - A(x) | : |x| \le 1 \}$$
such that the bound tends to $ 0$ when $K \to 1 .$

\begin{nonsec} \label{2.10} Canonical ring domains
\end{nonsec}
There are two ring domains in $\Bbb R^n$  
whose capacities are frequently used as
comparison functions. These canonical ring domains are the Gr\"otzsch ring
$R_G(s),\ s>1$, with complementary components $\overline B^n$ and
$\{te_1: t\ge s \}$ and the Teichm\"uller ring $R_T(t),\ t>0$, with 
complementary components $[-e_1,0]$ and $ [te_1,\infty)$. For the Gr\"otzsch
(Teichm\"uller) ring the capacity is the modulus of the curve family
joining the complementary components, denoted by $\gamma_n(s)$ and
$\tau_n(t),$ respectively. These capacities are related by
\begin{equation} \label{eq2.11}
\gamma_n(s)=2^{n-1}\tau_n(s^2-1)
\end{equation}
for $s>1$. There are several estimates for $\gamma_n(s)$ and
$\tau_n(t)$, for all $n\ge 3$; see \cite{G1}, \cite{A1}, \cite{AVV6},
\cite[ Section 7]{Vu1}.
When $n=2$ both functions can be expressed in terms of elliptic integrals;
see (\ref{eq2.23}) and (\ref{eq2.24}) below.

\begin{nonsec} \label{2.12} Conformal invariants $\mu_G$ and $\lambda_G$
\end{nonsec}
Since we are seeking
invariant formulations, the absolute
ratio is a natural tool. Another possibility is to use point-pair
invariants of a domain $G\subset \Bbb R^n,$ such as $\mu_G(a,b)$ or
$\lambda_G(a,b)$, $a,b\in G$, defined as follows
\begin{equation} \label{eq2.13}
\cases{\mu_G(a,b)=\inf_{C_{ab}}M( \Delta(C_{ab},\partial G;G)),\cr
       \lambda_G(a,b)=\inf_{C_a,C_b}M(\Delta(C_a,C_b;G)),\cr}
\end{equation}
where the infima are taken over all continua $C_{ab}$ (pairs of
continua $C_a,C_b)$) in $G$ joining $a$ and $b$ ($a$ to $\partial G$ and
$b$ to $\partial G$, resp.) Both $\mu_G$ and $\lambda_G$ are solutions
of the respective con\-for\-mal invariant extremal problems, and they both
have proved
to be efficient tools in the study of distortion theory of
quasi\-con\-for\-mal maps \cite{LF}, \cite{Vu1}.
Both $\mu_G$ and $\lambda_G^{-1/n}$ (\cite{LF})
are metrics for most subdomains of $\Bbb R^n$ (for $\mu_G$ we must
require that cap$\,\partial G>0$ and, for $\lambda_G,$ that card $(\overline
{\Bbb R}^n\setminus G)\ge 2)$.
\vskip 0.5cm
\centerline{\psfig{figure=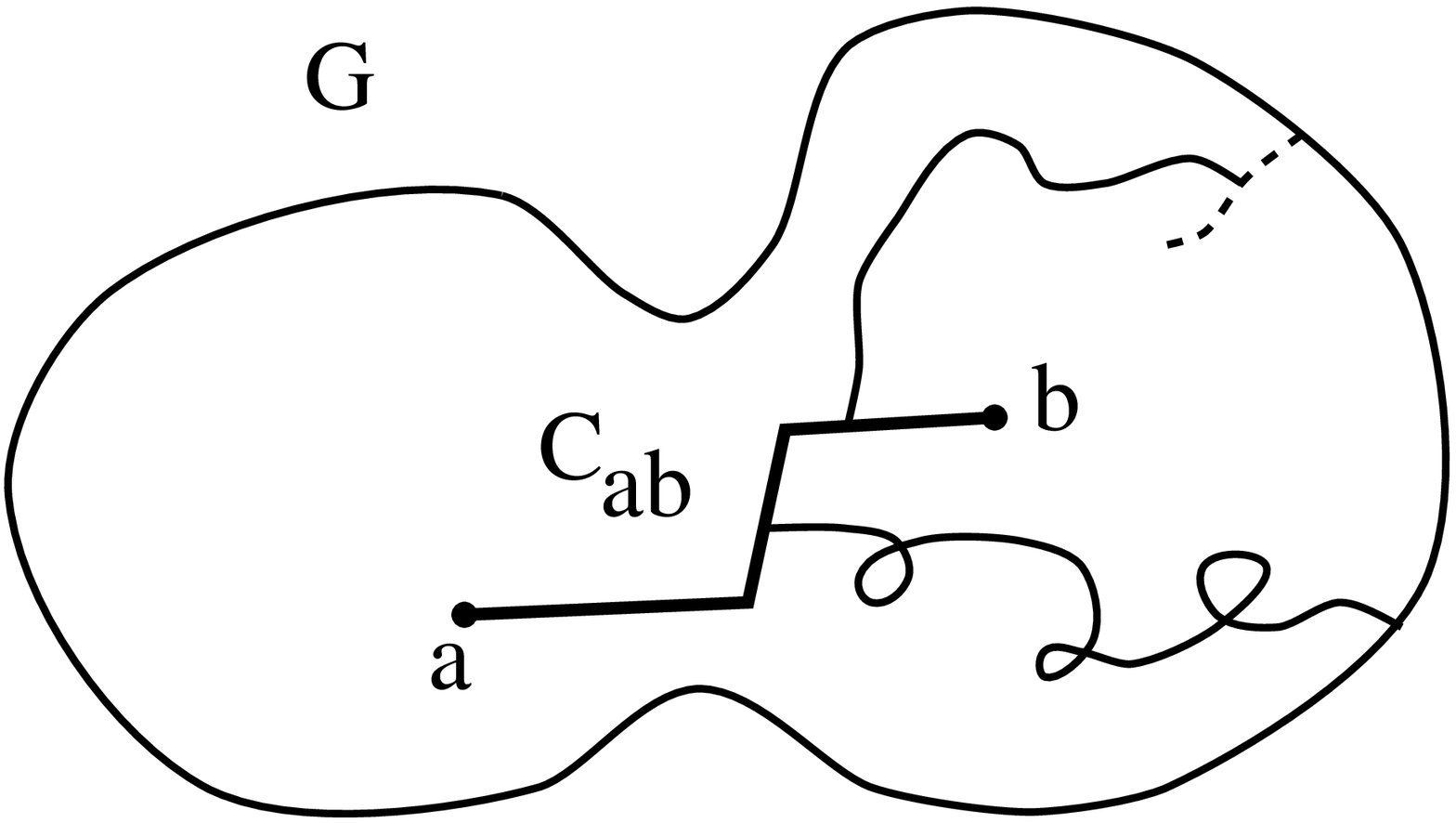,height=4cm}}
\vskip 0.5cm

\vskip 0.5cm
\centerline{\psfig{figure=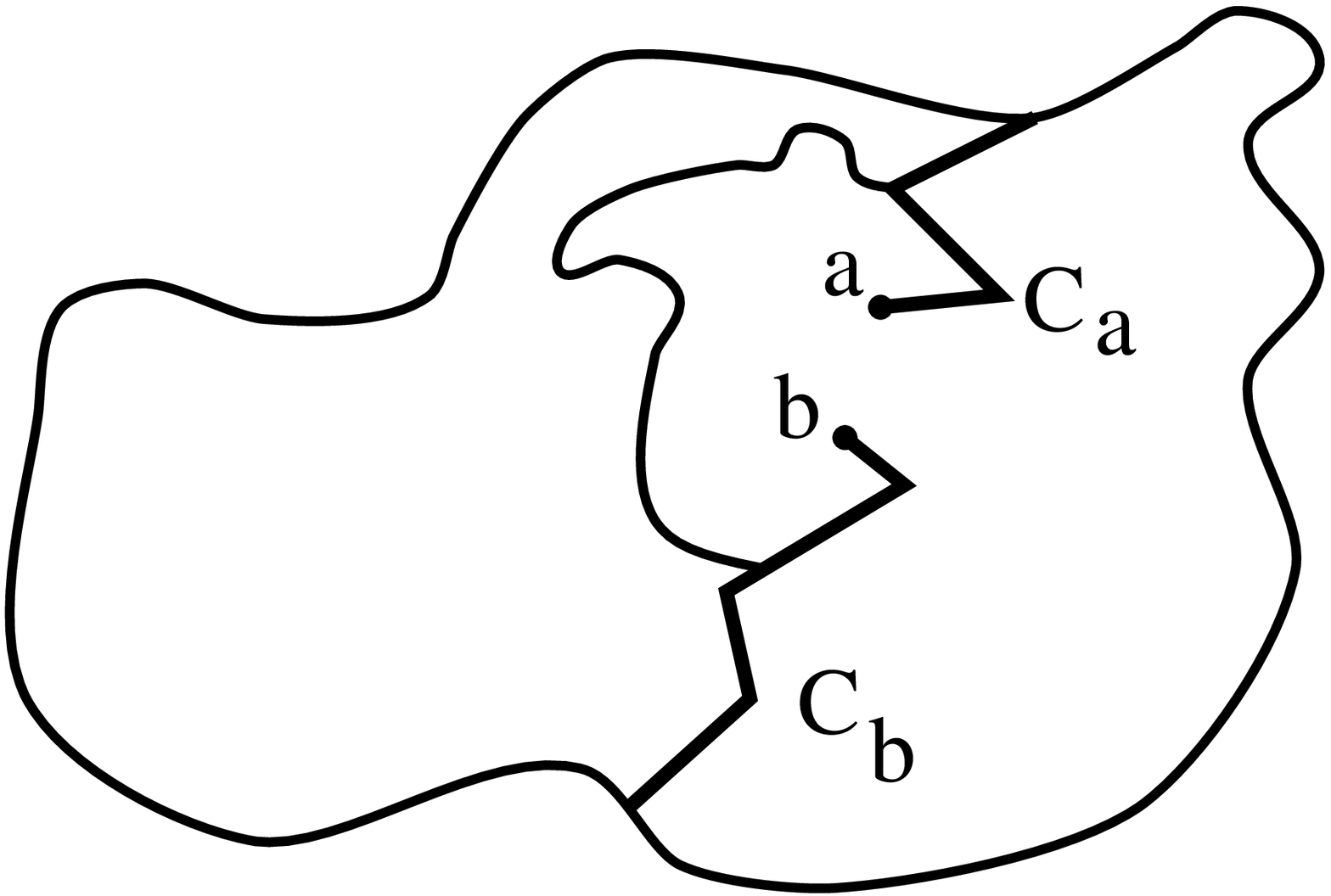,height=4cm}}
\vskip 0.5cm
\begin{nonsec} \label{2.14}
Bounds for $\mu_G$ and $\lambda_G$
\end{nonsec}
When comparing the mutual advantages of the absolute ratio and the extremal
quantities (\ref{eq2.13}) we note that the former is 
more explicit. On the other
hand---and this is the most important property of $\mu_G$ and 
$\lambda_G$---the transformation rules of $\mu_G$ and $\lambda_G$ under
$K$-quasi\-con\-for\-mal maps of the proper subdomain $G$ of $\overline{\Bbb
R}^n$ are simple: quasi\-con\-for\-mal maps are bilipschitz in the respective
metrics.
\begin{theo} { (Transformation rules)} \label{th2.15}
 If $f:G\to G'$ is $K$-quasi\-con\-for\-mal, then
$$(1)\ \mu_{fG}(f(a),f(b))/K\le \mu_G(a,b)\le K\mu_{fG}(f(a),f(b))\, ,$$
$$(2)\ \lambda_{fG}(f(a),f(b))/K\le \lambda_G(a,b)\le
K\lambda_{fG}(f(a),f(b))\, ,$$
for all $a,b\in G$.
\end{theo}

For all applications of these transformation rules we need estimates or
explicit formulas for $\mu_G$ and $\lambda_G$. Below we review what is
currently known about these invariants and point out some open problems.
Some applications of these invariants are given in \cite{F1} -\cite{F4},
\cite{Pa1}, \cite{Pa2}, \cite{Vu1}, \cite{AVV6}. 
Note that the important Schwarz lemma
for quasiconformal maps, Theorem \ref{th2.26},
 follows from Theorem \ref{th2.15}.

If $G=B^n$, then there is an explicit formula for $\mu_{B^n}$ as
well as one for $\lambda_{B^n}$ \cite{Vu1}.
For $n=2$, $G=\Bbb R^2\setminus \{0\}$, there is an explicit formula
for $\lambda_G$
which follows easily from the formula for the solution of the
Teichm\"uller extremal problem \cite[p. 192]{Kuz}.

Next, the general chart of inequalities among various metrics is given
in \cite{Vu0}.
If $f:\Bbb R^n\to \Bbb R^n$ is $K$-quasi\-con\-for\-mal and $G=fB^n\, ,$ then
there are upper and lower bounds for both $\lambda_G$ and $\mu_G$
\cite{Se1}. Next, if $G=\Bbb R^n\setminus \{0\}$, then there are
upper and lower bounds for $\lambda_G$ \cite[Section 8]{Vu1}.
If $\partial G$ is uniformly perfect, then there are lower bounds for
$\mu_G$ in terms of the metric $\delta_G$ in (\ref{eq2.6})---such 
a bound can be derived
from the results of \cite{JV}. Finally, for $G=B^2\setminus \{0\}$
there are upper and lower bounds for $\lambda_G$ \cite{LeVu}.

\begin{nonsec} \label{2.16} Open problem
\end{nonsec}
Find an explicit formula for $\lambda_{B^2\setminus \{0\}}$. Improve the
upper and lower bounds for $\lambda_G, \, G ={\Bbb R^n\setminus \{0\}}$.

\begin{nonsec} \label{2.17} Remark
\end{nonsec}
J. Ferrand proved in \cite{LF} that $\lambda_G^{-1/n}$ is a metric.
In \cite{AVV3} it
was shown that $\lambda_{B^n}(x,y)^{1/(1-n)}$ is a metric
and asked whether $\lambda_{B^n}(x,y)^{1/(1-n)}$ is a metric for more
general domains. Affirmative solutions were
subsequently found by A. Yu. Solynin \cite{Sol}, J. Jenkins \cite{J2},
and J. Ferrand
\cite{F3}.

\begin{nonsec} \label{2.18} Lipschitz conditions with respect to $\mu_G$
and $\lambda_G$
\end{nonsec}
The transformation rules of Theorem \ref{th2.15}
 are just special case of the
more general inequality (*). However, in many cases it is enough to use
Theorem \ref{th2.15} instead of (*). Therefore, the following question is
natural. Consider homeomorphisms $f:B^n\to f(B^n)=B^n$ satisfying the
property \ref{th2.15}(1) (or \ref{th2.15}(2)).
Are such maps quasi\-con\-for\-mal? This
question was raised by J. Ferrand \cite{LF} and a negative answer was
given in \cite{FMV}, where it was also shown that such maps are
H\"older-continuous.

\begin{nonsec} \label{2.19}
Heuristic principle
\end{nonsec}
The practitioners of quasi\-con\-for\-mal mapping theory have observed the
following heuristic principle: estimates for the modulus of a curve
family associated with a geometric configuration often lead to information
about quasi\-con\-for\-mal mappings. Unfortunately 
explicit formulas are available only in the simplest cases.
Symmetrization has proved
to be a very useful method for finding lower bounds for the solutions
of extremal problems such as the minimization of the capacities
of some suitable class of ring domains; see \cite{Ba2}, \cite{Ba3},
 \cite{Dub}, \cite{SolV}.
It should be noted that for dimensions $n\ge 3$ there is not even a
simple algoritm for the  numerical computation of the Gr\"otzsch capacity
$\gamma_n(s)$. For $n=3$ some computations were carried out
\cite{SamV}.
For the dimension $n=2$ there is an explicit formula for the Gr\"otzsch
capacity in terms of elliptic integrals, as we shall see below. We shall
also see that the same special functions will occur in several
function-theoretic extremal problems, and some most beautiful identities
for these special functions can be derived from Ramanujan's work on
modular equations.

In harmony with the above heuristic principle we now start a review of
special functions that will occupy a considerable part 
of Sections 2 and 3.

\begin{nonsec} \label{2.20}
Hypergeometric functions
\end{nonsec}
For $a,b,c\in \Bbb R,\ c\not= 0,-1,-2,$.. the (Gaussian) {\it
hypergeometric} {\it function} is defined by the series
$$F(a,b;c;r)=\sum_{n=0}^{\infty}{(a,n)(b,n)\over (c,n)n!}r^n$$
for $|r|<1,$ where $(a,0)=1,\ (a,n+1)=(a,n)(a+n),\ n=0,1,2,\dots$.
The hypergeometric function, one of the most important special
functions, was studied extensively by several eminent nineteenth
century mathematicians such as K.F. Gauss, E. Kummer, B. Riemann, H.A.
Schwarz, E. Goursat, and F. Klein \cite{Ask1}, \cite{Dut}, \cite{Kl2}.
Its importance is, in part, connected with its numerous particular
cases: there are lists in \cite{PBM} with hundreds of special cases of
$F(a,b;c;r)$ for rational triples $(a,b,c)$. Another reason for the
importance of $F(a,b;c;r)$ is its frequent occurrence in  several
different contexts in the 1990's, see \cite{Ao}, \cite{Ask2}, \cite{CC},
\cite{CH}, \cite{DM}, \cite{GKZ}, \cite{Var}, \cite{Va}, \cite{WZ1},
\cite{WZ2}. For our purposes, the main particular case of the
hypergeometric function is the {\it complete} {\it elliptic} {\it
integral} ${\K}(r)$ \cite{AS}, \cite{C3}, \cite{WW}
\begin{equation}
\label{eq2.21}
{\K}(r)={\pi\over 2} F({1\over 2},{1\over 2};1;r^2), \, 0\le r<1.
\end{equation}

\begin{nonsec} \label{2.22} Conformal map onto a disk minus a radial slit
\end{nonsec}
A con\-for\-mal mapping of a concentric annulus onto a disk minus a radial
segment starting from the origin is provided by an elliptic function.
The length of such a segment depends on the ratio of the radii in a
nonelementary fashion. In fact, if the inner and outer radius of the
annulus are $t\in (0,1)$ and 1, then the length $r\in (0,1)$ of the
radial segment satisfies the following transcendental equation, obtained
by equating the capacities of these two ring domains;
\begin{equation} \label{eq2.23}
{2\pi\over \log{1\over t}}={2\pi\over {\mu(r)}};\ \mu(r)={\pi\over 2}{
{\K}(r')\over {\K}(r)},
\end{equation}
where $r'=\sqrt{1-r^2}$ and we set $\mu(1)=0$. For $n=2$ the  Gr\"otzsch
capacity can be expressed as
\begin{equation} \label{eq2.24}
\gamma_2(s)=2\pi/\mu(1/s),\ s>1.
\end{equation}

\begin{nonsec} \label{2.25}
Schwarz lemma for quasi\-conformal maps
\end{nonsec}
The Schwarz lemma for analytic functions is one of the basic results of
complex analysis. A counterpart of this result also holds for
quasi\-con\-for\-mal maps in the following form.

\begin{theo} {}
\label{th2.26}
Let $f:B^n\to fB^n\subset B^n$ be $K$-quasi\-con\-for\-mal and $f(0)=0$. Then,
for $x\in B^n,$
$$(1)\ |f(x)|\le \varphi_{K,n}(|x|)\le
\lambda_n^{1-\alpha}|x|^{\alpha},\ \alpha=K^{1/(1-n)},$$
$$(2)\ |f(x)|\le
\psi_{K,n}(|x|)\equiv \sqrt{1-\varphi_{1/K,n}(\sqrt{1-|x|^2})^2},$$
where $\varphi_{K,n}(r) \equiv 1/\gamma_n^{-1}(K\gamma_n(1/r))$ and
$\varphi_{K,2}(r)=\mu^{-1}(\mu(r)/K)$. If, moreover, $fB^n=B^n$,
then
$$(3)\ |f(x)|\ge \varphi_{1/K,n}(|x|) \ge \lambda_n^{1- \beta} |x|^{\beta},$$
$$(4)\ |f(x)|\ge \psi_{1/K,n}(|x|).$$
\end{theo}

Note that in Theorem \ref{th2.26} both (1) and (2) are asymptotically sharp when
$K\to 1$. Here $\lambda_2 =4, \, \lambda_n\in[4,2e^{n-1})$ is a constant \cite{A2}. It
can be shown that, in (1) and (2),
$\varphi_{K,n}(r)$ and $\psi_{K,n}(r) $ are different for $n>2$ 
and identically equal for $n=2$.

\begin{nonsec} \label{2.27} Corollary \end{nonsec}
{\em If $f:B^n\to fB^n\subset B^n$ is $K$-quasi\-con\-for\-mal, then, for all
$a,b\in B^n,$
$$\tanh {\rho(f(a),f(b))\over 2}\le \varphi_{K,n}(\tanh{\rho(a,b)\over
2}).$$
}

It should be observed that in Corollary \ref{2.27} we do not require the
normalization $f(0)=0$. Corollary \ref{2.27} can be extended also to
domains $G$ of the form $f_1B^n$ where $f_1:\Bbb R^n\to \Bbb R^n$ is
$K_1$-quasi\-con\-for\-mal \cite{Se1}.

\begin{nonsec} \label{2.28}
Open problem
\end{nonsec}
Let $f:B^2\to B^2=fB^2$ be $K$-quasi\-con\-for\-mal, $a,b,c\in B^2$ and
$\alpha,\beta,\gamma$ the angles of the hyperbolic triangle with
vertices $a,b,c$. Find bounds for the hyperbolic area of the triangle
with vertices $f(a),f(b),f(c)$.

\begin{nonsec} \label{2.29}
Quasisymmetric maps $n=2$
\end{nonsec}
We recall that quasisymmetric maps already were defined and briefly
discussed in Section 1, where we pointed out that quasi\-con\-for\-mal maps of
$\Bbb R^n$ are $\eta$-quasisymmetric with an explicit $\eta_{K,n}$ given there.
For $n=2$ one can 
sharpen this result considerably, 
since a simple expression for $\eta_{K,2}$ is available by
the following result of S. Agard \cite{Ag}.

\begin{theo} {}
\label{th2.30}
A $K$-quasi\-con\-for\-mal map $f:\Bbb R^2\to \Bbb R^2$ is
$\eta_{K,2}$-quasisymmetric with
$$\eta_{K,2}(t)={u^2\over 1-u^2};\ 
u=\varphi_{K,2}\left(\sqrt{t\over 1+t}\right).$$
\end{theo}
Some sharp growth estimates for quasisymmetric maps were found by J.
Zaj\c{a}c in \cite{Za1}, \cite{Za2} in terms of the function
$\varphi_{K,2}(r)$. A related topic is the so-called Douady-Earle
extension problem, where sharp bounds were recently found by D. Partyka
\cite{Par3} in terms of the function $\varphi_{K,2}(r)$. The function
$\varphi_{K,2}$ satisfies many inequalities, which are sometimes used in
these studies. Some inequalities are given, e.g., in \cite{AVV2},
\cite{QVV2}, \cite{QVu1}.

\begin{nonsec} \label{2.30}
Schottky's theorem
\end{nonsec}
Schottky's classical result asserts the existence of a function
$\psi:(0,1)\times (0,\infty)\to (0,\infty)$ such that
$$\sup\{|f(z)| :\ |z|=r,\, f \in {\cal A}(t) \} \equiv \psi(r,t) \, ,$$
where
$${\cal A}(t)= \{ f: B^2 \to \Bbb R^2 \setminus \{0,1\}:  
f \mbox{ analytic,} \, |f(0)| =t \}.$$
Numerous explicit bounds for $\psi(r,t)$ have been found. W.K. Hayman
proved that
\begin{equation} \label{eq2.31}
\log\psi(r,t)\le (\pi+\log^+t) {1+r\over 1-r},
\end{equation}
and J. Hempel \cite{Hem1}, \cite{Hem2} proved, 
using some results of S. Agard \cite{Ag}, that
\begin{equation} \label{eq2.32}
\psi(r,t)=\eta_{M, 2}(t)\, , M = {1+r\over 1-r}.
\end{equation}
G. Martin \cite{Ma} found a new proof of (\ref{eq2.32}) based 
on holomorphic motions. One
can use Theorem \ref{th2.30} and (\ref{eq2.32})
to find sharper forms of Hayman's result
(\ref{eq2.31}) as shown in \cite{QVu2}. In \cite{QVu2} references to
related work by Jenkins, Lai, and Zhang are given.
Perhaps the best explicit estimate
known today is due to S.-L. Qiu \cite{Q3}: with 
$B=\exp(2\mu(1/\sqrt{1+t}))\, ,$
\begin{equation} \label{eq2.33}
16\eta_{K,2}(t)\le \min\{16t+B^K-B,\ (16t+8)^K-8\},
\end{equation}
for all $t\ge 0, \, K \ge 1$. Here equality holds if $K=1$ or $t=0$.

\begin{nonsec} \label{2.34}
Implementation of the heuristic principle
\end{nonsec}
We now give  an explicit example of the implementation of our heuristic
principle to the Schwarz lemma \ref{th2.26}. 
Indeed, the property that
$\mu(r)+\log r$ is monotone decreasing on $(0,1)$ implies the upper
bound $\varphi_{K,2}(r)\le 4^{1-1/K}r^{1/K}$ for $K>1$, $r\in (0,1)$
(this is the bound in Theorem \ref{th2.26}(1) with $n=2$ ).
Several other inequalities for $\varphi_{K,2}(r)$ can be proved in the
same way, if one uses other monotone functions involving $\mu(r)$ in
place of $\mu(r)+\log r$. A natural question is now: how do we find
such monotone functions? There is no simple answer to this question
beyond the obvious one: by studying ${\K}(r)$ and related functions
(recall that $\mu(r)= \pi{\K}(r')/(2 {\K}(r))$).
Many monotonicity properties of ${\K}(r)$ were found in the
1990's with help of ad hoc techniques from  classical analysis. What is
missing is a unified approach for proving such monotonicity results.
Since the publication of \cite{AVV2} in 1988 many inequalities and
properties for $\varphi_{K,n}(r)$ were obtained; see e.g.
\cite{QVV1}, \cite{QVV2}, \cite{QVu1}, \cite{QVu2}, \cite{Par1}-\cite{Par3},
\cite{Q3}, \cite{QV}, \cite{Za1}-\cite{Za2}. These papers solve many of
the open problems stated in \cite{AVV3}, \cite{AVV5}, and elsewhere. The
above heuristic principle has also found many applications there.
Another application of this principle occurs in \cite{AVV1}, where it was
shown that, for dimensions $n \ge 3 ,$ one can prove
many results for quasiconformal maps in a dimension-free way.
%Of the current open problems we mention here the
%following one from \cite{QVV2}.

\begin{nonsec} \label{2.35}
Open problem (\cite{QVV2})
\end{nonsec}
Show, for fixed $K>1, \, K \neq 2,$ that the function
$$g(K,r)={\mbox{artanh} \varphi_{K,2}(r)\over {\mbox{artanh}(r^{1/K})}}$$
is monotone from $(0,1)$ onto $(c(K),d(K))$
with $c(K) = \min \{K, 4^{1-1/K} \},
d(K) = \max \{K, 4^{1-1/K} \} .$ Note that $g(2,r) \equiv 2  $
since $\varphi_2(r) = 2 \sqrt{r}/(1+r) .$

\begin{nonsec} \label{Mori}
Mori's theorem
\end{nonsec}
A well-known theorem of A. Mori \cite{Mor1} states that a $K$-quasiconformal
map $f$ of the unit disk $B^2$ onto itself, normalized by $f(0) = 0,$ satisfies
\begin{equation} \label{eqMori}
 |f(x) -f(y)| \le M |x-y|^{1/K} \, ,
\end{equation}
for all $x,y \in B^2$
where $M =16 $ is the smallest constant independent of  $K.$
Clearly, this result is far from sharp if $K$ is close to $1. $
 O. Lehto and K. I. Virtanen \cite{LV} asked whether 
(\ref{eqMori}) holds with the constant $M = 16^{ 1-1/K} .$ This problem
has been studied by several authors, including 
R. Fehlmann and M. Vuorinen  \cite{FV} (the case $n \ge 2$),
V. I. Semenov \cite{Sem1},
S.-L. Qiu \cite{Q1}, G. D. Anderson, M. K. Vamanamurthy and M. Vuorinen
\cite{AVV4}. See also \cite{BP3}.
 Currently it is known that we can choose $M \le 64^{1 - 1/K} $
\cite[5.8]{AVV4}. Estimates of the function $\varphi_K(r)$ play a crucial
role in such studies.

\begin{nonsec} \label{MoriProb}
Open problem (from \cite{LV})
\end{nonsec}
Show that inequality (\ref{eqMori}) holds with $M \le 16^{1 - 1/K} .$
(Even the particular case when the points $x$ and $y$ are on the same
radius is open.)

\bigskip
Our next goal is to describe the classical method of computing
${\K}(r)$  in terms of the arithmetic-geometric mean. This procedure naturally
 brings forth the question of finding inequalities for ${\K}(r)$
in terms of mean values, which we will also touch upon.

\begin{nonsec} \label{2.36}
Arithmetic-geometric mean
\end{nonsec}
For $x,y>0$, the {\it arithmetic} and {\it geometric} {\it means} are
denoted by
$$A(x,y)=(x+y)/2,\ G(x,y)=\sqrt{xy},$$
respectively, and the {\it logarithmic} {\it mean} is defined by
$$L(x,y)={{x-y}\over {\log(x/y)}},\ x\not=y,\ L(x,x)=x.$$
Next, for $a>b>0$ let
$$a_0=a,\ b_0=b,\ a_{n+1}=A(a_n,b_n),\ b_{n+1}=G(a_n,b_n).$$
Then $a_n\ge a_{n+1}\ge b_{n+1}\ge b_n$ and
$$AG(a,b)\equiv \lim a_n=\lim b_n$$
is the  {\it arithmetic-geometric} {\it mean} of $a$ and $b$.
See \cite{BB1}, \cite{AlB}, \cite{ACJP}. Recently,
 the arithmetic-geometric mean has
been studied, in particular, in connection with the high-precision
computations of the decimal places of $\pi$ \cite{BBBP}, \cite{Lei}.
 The next theorem was proved by
Lagrange and Gauss (independently) some time between 1785 and 1799; see
\cite{C1} and \cite{Co}.
\begin{theo} {}
\label{th2.37}
For $r\in (0,1),\ r'=\sqrt{1-r^2}$,
we have
$${\K}(r)={\pi\over 2AG(1,r')}.$$
\end{theo}

For the approximation of ${\K}(r)$ in terms of mean values
it will be expedient to have notation as follows:
\begin {equation} \label{eq2.38}
M_t(x,y)=M(x^t,y^t)^{1/t},\ t>0,
\end{equation}
 for the modification of the means
$M=A,G,L,AG.$
These increase with $t$. There are numerous inequalities among the
above mean values, see
\cite{BB1}, \cite{BB2}, \cite{C2}, \cite{CV}, \cite{VV1}, \cite{San}.
The inequality $L(x,y)\le AG(x,y)$ for $x,y>0$ occurs in \cite{CV}.
In the opposite direction, the following theorem was proved by J. and
P. Borwein in 1994 \cite{BB2}.
\begin{theo} {}
\label{th2.39}
$AG(x,y)\le L_{3/2}(x,y)$ for $x,y>0$.
\end{theo}

%\begin{nonsec} \label{2.40} Open problem (from \cite{VV1})
%\end{nonsec}
%Is it true that $AG_t(x,y)\ge L(x,y)$ for all $x,y>0$ and all $t\in
%(t_0,1)$ for some $t_0\in (0,1)$?
%Can we choose $t_0=0.8$?

\bigskip
Some approximations for ${\K}(r)$ in terms of elementary functions
can be obtained if we use Theorem \ref{th2.39} or Theorem \ref{th2.37}
 and carry
out a few steps of the $AG$-iteration. In the next theorem two such
approximations are given. Part (1) is due R. K\"uhnau \cite{Kyuh3},  
and part
(2) is due to B.C. Carlson and J.L. Gustafson \cite{CG1}. See also
\cite{QV}.
\begin{theo} {}
\label{th2.41}
For $r\in[0,1)$ we have
$$(1)\ {\K}(r)>{9\over 8+r^2}\log{4\over r'},$$
$$(2)\ {\K}(r)<{4\over 3+r^2}\log{4\over r'}.$$
\end{theo}

For some recent inequalities for elliptic integrals, see \cite{AQV}.

%\begin{nonsec}{Five point problem} \label{2.45} \end{nonsec}
%For points $a,b,c,d \in {\Bbb R}^n$ recall that the absolute ratio is given by
%$$|a,b,c,d|=%{{q(a,c)q(b,d)}\over{q(a,b)q(c,d)}} =
%{|a-c||b-d|\over |a-b||c-d|}.$$
%In the limiting case when one of the points is $\infty$ we have
%$$|a,b,c,\infty| = {|a-c| \over |a-b|},\quad 
%|\infty,b,c,d| = {|b-d| \over |c-d|}. $$
%Hence:
%$$|a,b,c,d|= |a,b,c,\infty| |\infty,b,c,d|.$$
%\bigskip
%\begin{nonsec} 
%{Open problem} \label{2.46} \end{nonsec}
% Let $f : R^n \to R^n$ be $K$-quasicon\-for\-mal, and let
%$a, b, c, d \in R^n$ distinct. Find an upper bound for
%$$ |f(a),f(b),f(c),f(d)| $$
%which is explicit and asymptotically sharp as $K \to 1 .$

%\bigskip
%\noindent

%\noindent
%{\bf Idea }  We know \cite{Vu2} explicit $\eta_{K,n}(t),$
%$\eta_{1,n}(t) = t$ such that
%$$ {|f(a)-f(c)| \over |f(a)- f(b)| } \le \eta_{K,n}({|a-c| \over |a-b| })$$
%and 
%$$ {|f(b)-f(d)| \over |f(b)- f(c)| } \le  \eta_{K,n}({|b-d| \over |c-d| }).$$
%Hence
%$$ |f(a),f(b),f(c),f(d)| \le  \eta_{K,n}({|a-c| \over |a-b| })
% \eta_{K,n}({|b-d| \over |c-d| }).$$
%which has the desired properties. See also \cite{V4}.

%\bigskip
%%\noindent

%%\noindent
%{\bf Question} Can we simplify this bound? What about the
%case $n=2?$ (The bilip case is easy.)

\bigskip
\centerline{
\psfig{figure=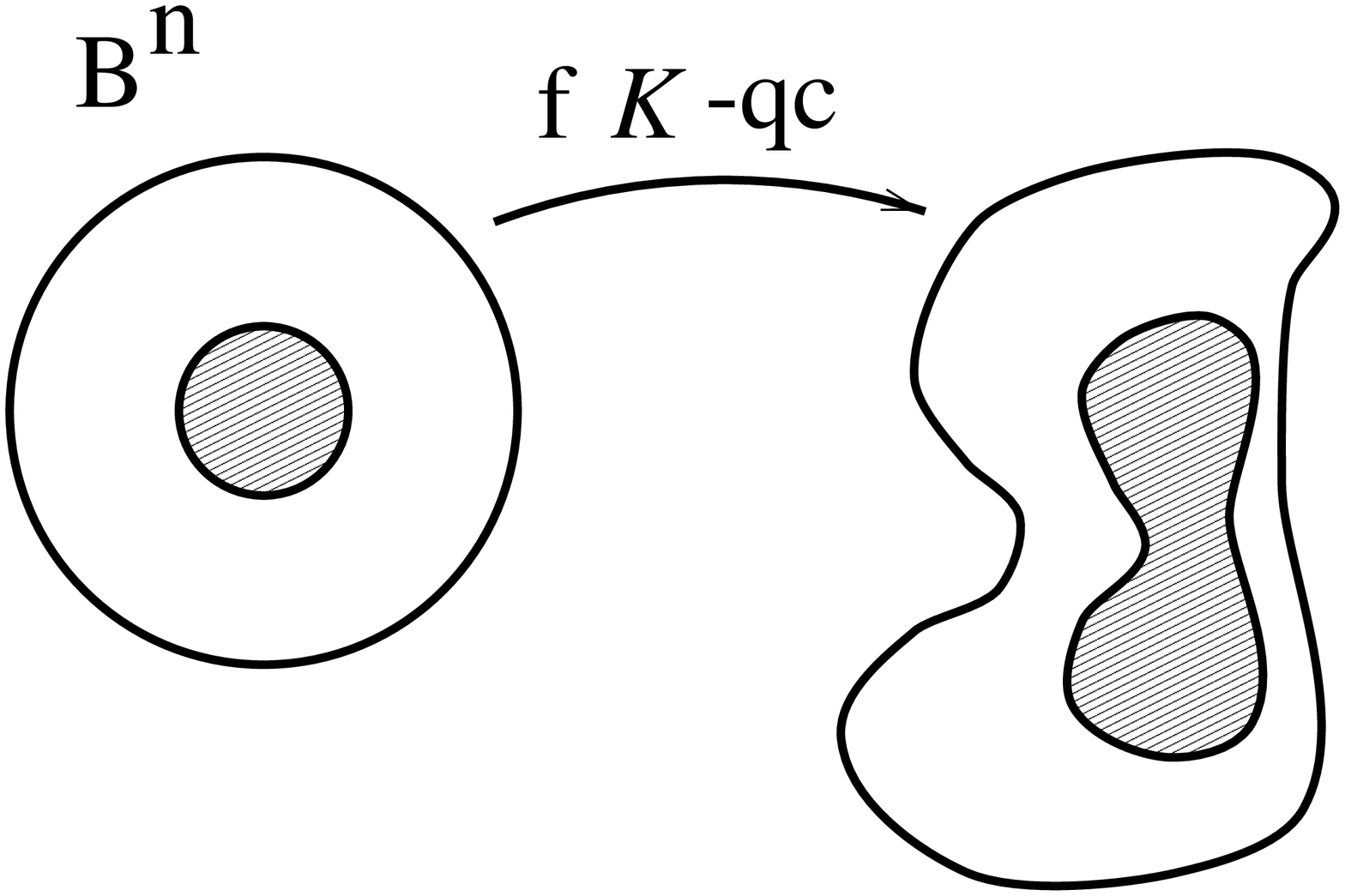,height=5cm}}

\bigskip
\begin{nonsec}
{Open problem ( from \cite[p. 193]{Vu1})} \label{2.50}\end{nonsec}
     Prove or disprove the following assertion. For each $n\ge 2$, $r\in
    (0,1)$, and $K\ge 1$ there exists a number $d(n,K,r)$ with $d(n,K,r)
    \to d(n,K)$ as $r\to 0$ and $d(n,K)\to 1$ as $K\to 1$ such that
    whenever $f : B^n \to {\Bbb R}^n$ is $ K$--qc, then $fB^n(r)$ is a
    ${d(n,K,r)}$--quasiball. More precisely, the representation
    $fB^n(r)=gB^n$ holds where $g: \overline{ \Bbb R}^n \to  
\overline{ \Bbb R}^n$ is a ${d(n,K,r)}$--qc
    mapping with $g(\infty)=\infty$. (Note: It was kindly pointed out by
    J.\ Becker that we can choose $d(2,1,r)=(1+r)/(1-r)$ either by
    \cite[pp.\ 39--40]{BC}  or by a more general result of S.\ L.\ Krushkal'
    \cite{KR}.)

\bigskip
%\noindent
Here we want an explicit constant $d(n,K)$---the existence follows e.g. from
the work of Tukia and V\"ais\"al\"a \cite{TV2}. An affirmative solution to this
local structure problem would have interesting applications.

%Fig. 28

%\eject

%%%%%%%%% SECTION 3 BEGIN
%%%%%%%%% SECTION 3 BEGIN
%%%%%%%%% SECTION 3 BEGIN
%%%%%%%%% SECTION 3 BEGIN
\cc
\section{  Recent results on special functions}
%MATERIAL FROM AVVB 
In this last section we shall mainly discuss recent results related to special
functions. We also mention a few geometric questions on quasi\-con\-for\-mal maps.

Some of the results below are related to the work of the Indian
mathematical genius S. Ramanujan   1887-1920. His published work has had
a deep impact on number theory, combinatorics, and special functions
(Hardy, Selberg, Dyson, Deligne). 

Ramanujan left numerous unpublished results in his notebooks at the time of
his premature death. It is estimated that the total number of his results
is 3000-4000. The publication of the edited notebooks with reconstructed
proofs by B. Berndt in 1985-1996 (5 volumes) made these results widely 
accessible. B. Berndt was awarded Steele Prize for Mathematical Exposition
in 1996 for this 
extraordinary achievement.

\begin{nonsec}
{\bf Asymptotic behavior of hypergeometric functions} \label{3.1}\end{nonsec}
The behavior of the hypergeometric function $F(a,b;c;r), \, a,b, c >0\, ,$
 at $r=1$ can be classified into three cases:

%\noindent
{\bf Case A.} $c> a+b.$ Now (Gauss) $ F(a,b;c;1) < \infty.$

%\noindent
{ \bf Case B.}  $c= a+b.  $ By Gauss' result as $r \to 1,$
$$ F(a,b;a+b;r) \sim {1\over B(a,b)} \log{1 \over 1-r} \,; B(a,b) = 
{\Gamma(a)\Gamma(b) \over \Gamma(a+b)} .$$

%\noindent
{\bf Case C.} $c< a+b. $ In this case the asymptotic relation is
$$ F(a,b;a+b;r) \sim D(1-r)^{c-a-b}, D = B(c,a+b-c)/B(a,b), $$
 as $r \to 1.$

The case $c=a+b$ is called the {\it zero-balanced case.}
%For rational triples $(a,b,c),$ $ F(a,b;c;r)$ sometimes reduces to a
%simpler function. Hundreds of special cases are given in \cite{PBM}.
The hypergeometric function satisfies numerous identities \cite{AS},
\cite{Bat1}-\cite{Bat3}. Perhaps the most famous are those due to
Kummer \cite{Pr}, \cite[15.5]{AS}.

\begin{nonsec} \label{3.2} {Ramanujan asymptotic formula}\end{nonsec}
The Gauss asymptotic formula in Case B was considerably refined by
Ramanujan, who proved \cite{Ask1} that if
 $\psi(a) = \Gamma'(a)/\Gamma(a)$ and
$$ R(a,b) =- \psi(a)-\psi(b) -2 \gamma, \quad R(1/2,1/2) = \log 16 ,$$
where
$\gamma$ stands for the Euler-Mascheroni constant and $B(a,b)$ for the
beta function, then
$$ B(a,b)F(a,b;a+b;r)+\log(1-r) = $$
$$R(a,b) + O((1-r) \log(1-r)) $$
as $r \to 1. $ 
%with $\psi(a) = \Gamma'(a)/\Gamma(a)$
%$$ R(a,b) =- \psi(a)-\psi(b) -2 \gamma, \quad R(1/2,1/2) = \log 16 ,$$
 
Ramanujan's result was extended recently, in \cite{ABRVV}
where it was shown e.g. that for $a, b  \in (0,1), \, B= B(a,b),$
$$B F(a,b;a+b;r) + (1/r) \log(1-r) $$
is increasing on $ (0,1)$ with range $(B-1,R).$
See also \cite{PV1} where $a+b$ is replaced by $c$. Convexity properties
of the hypergeometric function of the unit disk have been studied
recently e.g. in \cite{PV2}, \cite{PV3}, \cite{PSa}.

\begin{nonsec}{ Perturbation of identity} \label{3.3}\end{nonsec}
By continuity, small
changes of argument lead to small changes
of the values of the function.
For example,
we expect that
$$ F(a_1,b_1;c_1;r) \mbox{ and } F(a_2,b_2;c_2;r) $$
are close if the parameters $(a_1,b_1,c_1) \mbox{ and } 
(a_2,b_2,c_2)$ are close. In view of the asymptotic behavior at
$r=1$ considered in \ref{3.1} above, it seems natural to require
that $c_1-a_1-b_1 = c_2-a_2-b_2.$
Two natural questions are:
When are $\K(r)$ and $(\pi/2)F(a,b;a+b;r^2)$
close to each other (recall that $\K(r) = (\pi /2) F(1/2,1/2;1;r^2)$)? 
Can we extend the many properties of $\K(r)$ to
$F(a,b;a+b;r^2)?$

\begin{nonsec}{ Landen inequality} \label{3.4}\end{nonsec}
Recall first that one of the most important properties of $\K(r)$
is given by the {\em Landen identity} (1771) (see \cite{AlB}, \cite{Hou}),
which states that for all $r \in (0,1)$
\begin{equation} \label{3.5}
 \K({2 \sqrt{r}\over 1+r}) = (1+r)\K(r).
\end{equation}
The next theorem, an extension of (\ref{3.5}), might be called 
a {\em Landen inequality} \cite{QVu3}.

\begin{theo} {} \label{3.6}
 For $a,b \in (0,1)$ with
$a+b \le 1$ we have, for all $r \in (0,1),$
$$ F(a,b;a+b;({2 \sqrt{r}\over 1+r})^2) \le (1+r) F(a,b;a+b;r^2).$$
\end{theo}
%According to Askey, this is a new result.

In \cite{AVV5} there is a list of monotonicity properties of $\K(r).$
For instance, $\K(r)/\log(4/r')$ is monotone decreasing. It is
natural to ask if such properties have a generalization for hypergeometric
functions. The next open problem is of this type.

\begin{nonsec} {Open problem} \label{3.7}
\end{nonsec} Let $a,b \in (0,1]$ with $a+b \le 1, $
$$Q(r)= B(a,b)F(a,b;a+b;r)/\log{c \over 1-r}\, , \quad c= e^{R(a,b)}, $$
and $G(r)=(Q(r)-1)/(1-r).$ Is it true that the Maclaurin coefficients of
$G(r)$ are positive?

\bigskip
An affirmative answer would give a refinement of Ramanujan's asymptotic
formula \ref{3.2} and also imply that $G$ is strictly increasing and convex.
Computer experiments suggest that the answer is in the affirmative.

The last topic of this section deals with the algebraic
identities for the function
$$\varphi_K(r)=\varphi_{K,2}(r)=\mu^{-1}(\mu(r)/K)$$
that follow from Ramanujan's work on modular equations \cite{Bern1},
\cite{Bern2}, \cite{Bern3}. Of these \cite[pp. 8-9]{Bern3} contains a very
helpful list of Ramanujan's numerous contributions in the field. Because
Ramanujan's work in this field became widely accessible only with the
publication of \cite{Bern3} in 1991, the derivation of these results as
corollaries to Ramanujan's work could not have been possible before
1991.

\begin{nonsec} \label{2.42}
Modular equations of degree p
\end{nonsec}
The argument $r\in(0,1)$ of the complete elliptic integral ${\K}(r)$
is sometimes called the {\it modulus} of ${\K}$. A {\it modular}
{\it equation} of {\it degree} $p>0$ is the relation
\begin{equation} \label{eq2.43}
{{\K}(s')\over {\K(s)}}=p{{\K}(r')\over {\K(r)}}
\Leftrightarrow \mu(s)=p\mu(r).
\end{equation}
The solution of this equation is $s=\varphi_{1/p}(r)$.
Modular equations were studied by several mathematicians in the
nineteenth century. The most remarkable progress was made, however, by
Ramanujan in 1900-1920. We first record a few basic properties
of $\varphi_K(r)$ which will be handy for the discussion of
modular equations:
\begin{equation} \label{eq2.44}
\cases{\varphi_K(r)^2+\varphi_{1/K}(r')^2=1,\cr
       \varphi_A(\varphi_B(r))=\varphi_{AB}(r),\ A,B>0,\cr
        \varphi_{1/K}(r)=\varphi_K^{-1}(r),\cr
        \varphi_{2}(r)= \displaystyle{{2 \sqrt{r} \over 1+r }} , \cr }
\end{equation}
for all $r\in [0,1]$. The classical {\it Legendre-Jacobi modular equation of
order 3}
\begin{equation} \label{eq2.45}
\sqrt{rs}+\sqrt{r's'} =1,\ s=\varphi_{1/3}(r),
\end{equation}
can be solved for $s.$ We can easily find the solution if we use a
symbolic computation program such as Mathematica. The solution
was worked out (by hand!) in \cite{KZ}.

\begin{nonsec} \label{2.45b}
Ramanujan modular equations
\end{nonsec}
We use the term modular equation not only for the transcendental
equation (\ref{eq2.43}) but also for an algebraic equation that follows from
(\ref{eq2.43}), as in \cite{Bern3}.
An example of such an algebraic equation is the Legendre-Jacobi
modular equation (\ref{eq2.45}). We now rewrite (\ref{eq2.45}) using
Ramanujan's notation:
$$ \sqrt[4]{ \alpha  \beta } +\sqrt[4]{(1- \alpha)(1-  \beta) } =1,
 \quad \alpha = r^2, \, \beta = \varphi_{1/3}(r)^2,$$
for all $r \in (0,1) . $
Following \cite{Bern3} and \cite{Vu5}
we now give a few of Ramanujan's modular equations.

\begin{theo} {} \label{thRnew}  The function $\varphi_K$ satisfies the 
following identities:
\smallskip
\item{(1)} For $\alpha =r^2,\ \beta =\varphi_{1/5}(r)^2$, we have
$$(\alpha\beta )^{1/2}+\{ (1-\alpha )(1-\beta )\}^{1/2}+2\{
16\alpha\beta
(1-\alpha )(1-\beta )\}^{1/6}=1 . $$

\item{(2)} For $\alpha =r^2,\ \beta =\varphi_{1/7}(r)^2$, we have
$$(\alpha\beta )^{1/8}+\{ (1-\alpha )(1-\beta )\}^{1/8}=1 . $$

\item{(3)} For $\alpha =r^2,\ \beta =\varphi_{1/3}(r)^2,\ \gamma
=\varphi_{1/9}(r)^2$ we have $$\{\alpha (1-\gamma )\}^{1/8}+\{\gamma
(1-\alpha )\}^{1/8}=2^{1/3}\{\beta (1-\beta )\}^{1/24} .$$

\item{(4)} For $\alpha =r^2,\ \beta =\varphi_{1/23}(r)^2$, we have
$$(\alpha\beta )^{1/8}+\{ (1-\alpha )(1-\beta
)\}^{1/8}+2^{2/3}\{\alpha\beta
 (1-\alpha )(1-\beta )\}^{1/24}=1 .$$

\item{(5)} For $\alpha =r^2,\ \beta =\varphi_{1/7}(r)^2,\
$ or for $\alpha =\varphi_{1/3}(r)^2,\ \beta=\varphi_{1/5}(r)^2$, we
have
$$(\alpha\beta )^{1/8}+\{ (1-\alpha )(1-\beta )\}^{1/8}-\{\alpha\beta
(1-\alpha )(1-\beta )\}^{1/8}= $$
$$\{ {{\displaystyle 1}\over { \displaystyle
2}}(1+\sqrt{\alpha\beta}+\sqrt{(1-\alpha )(1-\beta )})\}^{1/2} .$$

\bigskip
\end{theo}

{\bf Proof.} All of these identities are from \cite{Bern2}: (1) is 
\cite[ p.280, Entry 13 (i)]{Bern2}; (2) is p. 314, Entry 19 (i); 
(3) is p. 352, Entry 3
(vi); (4) is p. 411, Entry 15 (i); and (5) is p. 435, Entry 21 (i). 

\bigskip

\begin{theo} {}
\label{th2.47}
 The function $\varphi_K$ satisfies the following
identities for $s\in (0,1):$
\smallskip

\indent (1) $xy+x'y'+2^{5/3}\{ xyx'y'\}^{1/3}=1$, %\hfill\break
where $x=\varphi_{\sqrt 5}(s),$ $\, y=\varphi_{1/\sqrt 5}(s)$,
\smallskip

\indent (2) $(xy)^{1/4}+(x'y')^{1/4}=1$, %\hfill\break
where $x=\varphi_{\sqrt 7}(s),$ $\, y=\varphi_{1/\sqrt 7}(s)$,
\smallskip

\indent (3)
$(xy)^{1/4}+(x'y')^{1/4}=2^{1/3}(s^2(1-s^2))^{1/24}$ %\hfill\break
where $x=\varphi_3(s),$ $\, y=\varphi_3(s')$,
\smallskip

\indent (4)
$(xy)^{1/4}+(x'y')^{1/4}+2^{2/3}(xx'yy')^{1/12}=1$ %\hfill\break
where $x=\varphi_{1/\sqrt{23}}(s),\,$ $ y=\varphi_{\sqrt{23}}(s')$,
\smallskip

\indent (5) $(xy)^{1/4}+(x'y')^{1/4}-\{ xx'yy'\}^{1/4}=$
\indent $\{ {{ 1}\over
{\ 2}}(1+xy+{x'y'})\}^{1/2}$ %\hfill\break
where $x=\varphi_{\sqrt{5/3}}(s),$ $\, y=\varphi_{\sqrt{3/5}}(s)$.
\end{theo}
%\end{nonsec}
\bigskip

{\bf Proof.} All parts follow from Theorem  \ref{thRnew} above in the
same way after $r$ is chosen appropriately. For this reason we give here
the details only for (5). For (5) set $r=\varphi_{\sqrt{15}}(s)$. By
(2) and (3) we see that
$$\alpha =\varphi_{1/3}(r)^2=\varphi_{\sqrt{5/3}}(s)^2,\quad 1-\alpha
=\varphi_{\sqrt{3/5}}(s')^2 \, ,$$
and thus the proof follows from Theorem \ref{thRnew} (5).

%\eject
\bigskip

\begin{nonsec} \label{2.48} {Corollary}\end{nonsec}
{\em  We have the following identities: }
$$2uu'+2^{5/3}(uu')^{2/3}=1;\quad u=\varphi_{\sqrt{5}}(1/\sqrt{2}),\leqno (1)$$
$$2(uu')^{1/4}=1;\quad u=\varphi_{\sqrt 7}(1/\sqrt{2}),\leqno (2)$$
$$\sqrt u+\sqrt{u'}=2^{1/4};\quad u=\varphi_3(1/\sqrt2),\leqno (3)$$
$$2(uu')^{1/4}+2^{2/3}(uu')^{1/6}=1;\quad u=\varphi_{\sqrt{23}}(1/\sqrt
2),\leqno (4)$$
$$2(uu')^{1/4}-(u u')^{1/2}=\{ {{1}\over {2}}(1+2uu')\}^{1/2};\leqno (5)$$
$$u=\varphi_{\sqrt{5/3}}(1/\sqrt 2).$$

\bigskip

{\bf Proof.} All parts follow from Theorem \ref{th2.47} and 
(\ref{eq2.44}) in the
same way. We give here the details only for (5). Set $s=1/\sqrt 2$ in
Theorem \ref{th2.47} (5) 
and observe that then, $x'=y,\ y'=x$ and thus (5) follows as
desired.

\bigskip

%{\bf 7. Notation.} For $|c|\in (0,1/\sqrt 2)$ let
%$$r(1,c)=\sqrt{(1+\sqrt{1-4c^4)/2}}$$
%$$r(2,c)=\sqrt{(1-\sqrt{1-4c^4)/2}}$$
%Then $u^2(1-u^2)=c^2$ with $u=r(j,c),\ j=1,2,$ and the upper bound of
%$|c|$ implies that $r(j,c),\ j=1,2,$ are real.
%\eject

%Parts (1)-(2) of Theorem \ref{th2.49} are given in [BB, p.139]
%whereas parts (3), (4) are perhaps new.

%\bigskip
%\begin{theo} {} \label{th2.49}
%The following identities hold:
%$$\varphi_{1/\sqrt 3}(1/\sqrt 2)={\sqrt{2-\sqrt 3}\over
%2}\leqno (1)$$
%$$\varphi_{1/\sqrt 7}(1/\sqrt 2)={\sqrt{8-\sqrt{63}}\over
%4}\leqno (2)$$
%$$\varphi_{1/\sqrt{5/3}}(1/\sqrt{2})=\sqrt{{1 \over 2} \left(1+
%\sqrt{ 1-4\left({1+ \sqrt{5} \over 4}\right)^8}\right)}\leqno (3)$$
%$$\varphi_{1/\sqrt{23}}(1/\sqrt{ 2})=\sqrt{{1-\sqrt{1-4z^{24}}\over
%2}}\leqno (4)$$
%where $z=-{{1}\over {3\root 3\of 2}}+{{1}\over {
%3(25+\sqrt{621})^{1/3}}}+{{ \root 3 \of {25+\sqrt{621}}}\over {
%3\cdot 2^{2/3}}}$.
%\end{theo}
%\bigskip

%{\bf Proof.} (1) Jacobi's classical modular equation gives
%$$\sqrt{r\varphi_{1/3}(r)}+\sqrt{r'\varphi_3(r')}=1$$
%and then set $r=\varphi_{\sqrt 3}(1/\sqrt 2)$. 
%Solving this equation yields the desired formula (1), and then (2)
%and (3) follows from Cor. 6 (2) and (4), respectively.
%For the last formula the Mathematica program was used.

\begin{nonsec} { Generalized modular equations} \label{3.8}\end{nonsec}
%Recall first that
%the classical modular equations 
%$$ {{\K}(s') \over {\K}(s)}= p {{\K}(r') \over {\K}(r)} 
%\Longleftrightarrow \mu(s) = p \mu(r) $$
%were considered in (\ref{2.43}).
A {\it generalized modular equation with signature }$1/a$ {\it and order (or
degree)} $p$ is
\begin{equation} \label{GenModEq}
 {F(a, 1-a;1;1-s^2)\over F(a, 1-a;1;s^2) } = 
p {F(a, 1-a;1;1-r^2)\over F(a, 1-a;1;r^2) }.
\end{equation}
Such equations were studied extensively by Ramanujan, who also gave a great
number of algebraic identities for the solutions.
Many of his results were proved in 1995 by Berndt, Bhargava, and Garvan
in a long paper \cite{BBG} (see also \cite{Gar}). 
The main cases they studied are:

\bigskip
{{ $$ a= 1/6, 1/4, 1/3, \quad p = 2,3, 5, 7, 11, ...$$ }  }
\bigskip

With
$$ \mu_a(r)={\pi \over 2 \sin (\pi a)}
{F(a, 1-a;1;1-r^2)\over F(a, 1-a;1;r^2)}$$
the solution of (\ref{GenModEq}) is given by
$$ s = \mu_a^{-1}(p \mu_a(r) ) \equiv \varphi^a_{1/p}(r).$$
Note that $\mu_a(r) = \mu_{1-a}(r)$ for $a \in (0,1/2) $ and
$ \mu(r) = \mu_{1/2}(r) .$

\bigskip
%Some examples of the results of the 70 pages long paper \cite{BBG} are next
%given.
For generalized modular equations the Ramanujan notation is

\bigskip
\centerline{
{ $\alpha \equiv r^2 , \quad  \beta \equiv \varphi^a_{1/p}(r)^2 $ } . }
\bigskip

\begin{theo} {\cite[Theorem 7.1]{BBG}}
\label{th7.1}
If $\beta$ has degree 2 in the theory of signature 3, then, with
 $a=1/3, \alpha= r^2, \beta= \varphi_{1/2}^a(r)^2\, , $
$$(\alpha\beta)^{1\over 3}+\{(1-\alpha)(1-\beta)\}^{1\over 3}=1.$$
\end{theo}

\begin{theo} {\cite[Theorem 7.6]{BBG}}
\label{th7.6}
If $\beta$ has degree 5 then,
with $a=1/3, \alpha= r^2, \beta= \varphi_{1/5}^a(r)^2\, , $
\end{theo}
\begin{equation} \label{eq7.17}
(\alpha\beta)^{1\over 3}+\{(1-\alpha)(1-\beta)\}^{1\over
3}+3\{\alpha\beta(1-\alpha)(1-\beta)\}^{1\over 6}=1.
\end{equation}

%\begin{theo} {}
%\label{th7.7}
%If $\beta$ has degree 7, then

%\end{theo}
%\begin{equation} \label{eq7.24}
%m=({\beta\over \alpha})^{1\over 3}+({1-\beta \over 1-\alpha})^{1\over 3}-
%{7\over m}({\beta(1-\beta)\over \alpha(1-\alpha)})^{1\over
%3}-3({\beta(1-\beta)\over \alpha(1-\alpha)})^{1\over 3}.
%\end{equation}

\begin{theo} {\cite[Theorem 7.8]{BBG}}
\label{th7.8}
If $\beta$ has degree 11 then,
with $a=1/3, \alpha= r^2, \beta= \varphi_{1/11}^a(r)^2\, , $
\end{theo}
\begin{equation} \label{eq7.28}
(\alpha\beta)^{1\over 3}+\{(1-\alpha)(1-\beta)\}^{1\over
3}+6\{\alpha\beta(1-\alpha)(1-\beta)\}^{1\over 6}+
\end{equation}
$$3\sqrt 3
\{\alpha\beta(1-\alpha)(1-\beta)\}^{1\over 12}\{(\alpha\beta)^{1\over
6}+\{(1-\alpha)(1-\beta)\}^{1\over 6}\}=1.$$

%Fig. 26

Several open problems are now immediate.

\begin{nonsec}{ Open problem} \label{3.41}\end{nonsec}
To what extent can the properties of $\mu(r)$
be extended for $\mu_a(r)$?

\begin{nonsec}{ Open problem} \label{3.42}\end{nonsec}
 To what extent can the properties of $\varphi_K(r)$
 be extended for $\varphi^a_K(r)$?

\bigskip

Solving these problems will require very extensive studies.
A basic tool is the {\em Ramanujan derivative formula} 
\cite[p. 86]{Bern2}:
\begin{equation} \label{RamDerF}
 {d \mu_a(r) \over dr } = -{1 \over r(1-r^2)} { 1 \over
F(a,1-a;1;r^2)^2  }. \end{equation}
A direct application of the $F(a,b;c;r)$ derivative formula 
$$ {d \over dr}F(a,b;c;r) =  {ab\over c} F(a+1,b+1;c+1;r)$$ 
from \cite[15.2.1]{AS}
leads to a more
complicated form than (\ref{RamDerF})
so the formula (\ref{RamDerF})
is specific for the case $b=1-a, c=1.$ By equating this more complicated
form and (\ref{RamDerF}) we obtain the following interesting identity
for $a,r \in (0,1):$
\begin{equation} \label{RamId}
\cases{F(1+a,2-a;2;1-r)F(a,1-a;1;r)+\cr
F(1+a,2-a;2;r)F(a,1-a;1;1-r)={\displaystyle
{\sin(\pi a)\over \pi a(1-a)r(1-r)}}. \cr}
\end{equation}

\begin{theo} {\cite{BPV}} \label{3.43}
For $0 <a \le 1/2, r ,s \in (0,1),$ we have
$$ \mu_a(r) + \mu_a(s) \le 2\mu_a(\sqrt{{2rs \over 1+rs + r' s'}} )
%{ 2 \sqrt{rs} \over \sqrt{(1+r)(1+s)} + \sqrt{(1-r)(1-s)}}) $$
\le 2\mu_a(\sqrt{rs}), $$
with equality for $r =s .$
\end{theo}
\bigskip
%\noindent
The above inequality $ \mu_a(r)+\mu_a(s) \le 2 \mu_a(\sqrt{rs})$
resembles the multiplicative property of the logarithm
$ \log a + \log b = 2 \log \sqrt{a b}\, , a, b > 0 ,$ and hence
$\mu_a(r)$ behaves, to some extent, like a logarithm.

\bigskip
%Note that
%$${ 2 \sqrt{rs} \over \sqrt{(1+r)(1+s)} + 
%\sqrt{(1-r)(1-s)}}= \sqrt{{2rs \over 1+rs + r' s'}} .$$

\bigskip
%\noindent
%{\bf Thm [QVu4]} 

\begin{theo} {\cite{QVu4}} \label{3.45} 
For $a \in (0, 1/2]$ and $ r, t \in (0,1),$
$$  2 \mu_a({r+t \over 1 + rt + r' t'}) \le  \mu_a(r) + \mu_a(t)\, , $$
with equality for $ t=r .$
\end{theo}

\bigskip

\begin{nonsec}{Open problem} \label{3.46}\end{nonsec}
In view of these results it is natural to ask if
there is an addition formula for $\mu_a .$
\bigskip

\begin{nonsec}{ Duplication inequality} \label{3.47}\end{nonsec}
The Landen identity yields the following duplication formula for $\mu(r)$
$$ \mu(r) =2 \mu({2 \sqrt{r} \over 1+r}) .$$
The next theorem from \cite{QVu4} 
could be called a {\em duplication inequality} for $\mu_a(r) .$
\bigskip
%\noindent
\begin{theo} {} \label{3.48}
 For $a \in (0,1/2]$ let $R = R(a, 1-a),$
$$ C\equiv \left( 1+{ \sin \pi a \over \pi}(R- \log 16) \right)^2 \, ,$$
and $C_1 = \min \{ 2, C \}. $ 
 Then, for $r \in (0,1),$
$$ \mu_a(r) \le 2\mu_a({2 \sqrt{r} \over 1+r}) \le C_1 \mu_a(r).$$
\end{theo}

\bigskip
%For the proof of these results several auxiliary results were proved.
%It seems likely that in near future there will be new results for the
%function $\varphi^a_K(r).$ At the moment there are hardly any
%inequalities for this function.

Jacobi's work yields dozens of infinite product expansions for elliptic
functions. A representative identity is the following one, where
$ q = \exp(-2 \mu(r) ) $
$$ \exp (\mu(r) + \log r) = 4 \prod_{n=1}^{\infty}\left( {1+q^{2n}\over
1+q^{2n-1}}\right)^4 .$$
Observe that $\mu(r)$ occurs on both sides!

\bigskip
%\noindent
\begin{theo} {} \label{3.50}
 For $a \in (0,1/2]$ let $R = R(a, 1-a).$ For $r \in
(0,1)$
let $$r_0 = \sqrt{1-r^2}, r_n = {2 \sqrt{r_{n-1}} \over 1+ r_{n-1}},
\quad p=\prod_{n=0}^{\infty} \left(1 + r_n \right)^{2^{-n}}.$$
Then
$$ p \le \exp( \mu_a(r) + \log r) \le {\exp R \over 16}p,
 $$
with equality for $a= 1/2 .$
\end{theo}
\bigskip

\bigskip
%\noindent
In this theorem from \cite{QVu5} 
even the case $a= 1/2$ is new. Note that $\mu_a(r)$
occurs only in the middle term!

\bigskip
%\noindent

\begin{nonsec}{ Linearization for $\varphi_K(r)$ } \label{3.51} 
\end{nonsec} It was observed in [SamV]
that the functions $p(x) =  \log(x/(1-x)), q(x) =e^x/(1+e^x)$,
$$p : (0,1) \to \Bbb R, \quad q: \Bbb R \to (0,1)\quad p = q^{-1} $$
have a regulating effect on $\varphi_K(r):$ 

\centerline{
\psfig{figure=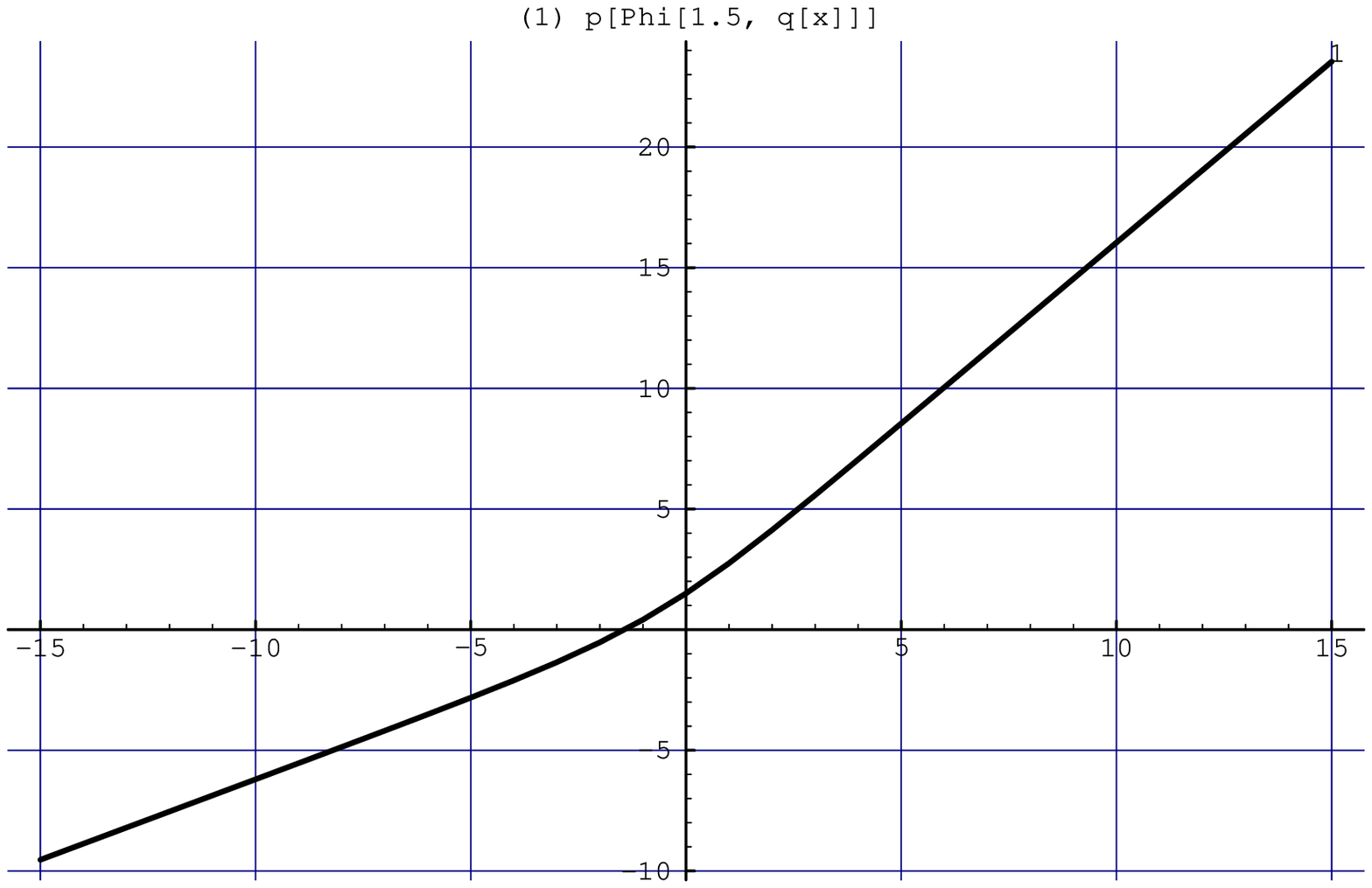,width=11cm}}

\bigskip
%\noindent
\begin{theo} {(\cite{AsVV})}
\label{3.53}The function $g : \Bbb R \to \Bbb R,
g(x)= p (\varphi_K(q(x))$ has 
increasing derivative with range $(1/K,K).$
\end{theo}

\bigskip
%\noindent
We call such a transformation {\sl a linearization}. 

%\centerline{
%\psfig{figure=lubpic27.ps,width=11cm}}

\bigskip
%\noindent
\begin{nonsec}
 { Open problem } \label{3.54}\end{nonsec}
Find a similar result for $\varphi^a_K(r). $

\bigskip
%\noindent

\bigskip
\begin{nonsec}
{Conformal Invariants Software } \label{3.62}\end{nonsec}
We shall discuss here briefly
the C, Mathematica, and MATLAB language software supplement for
the monograph \cite{AVV6}, which will soon be completed. This software
will fill one disk (1.4M), and it provides, in these languages, algorithms
for the special functions mentioned earlier. 
A useful survey of special function computation is \cite{LOl}.
 See also \cite{AS}, \cite{Bak},
\cite{Mosh}.

Nine example programs are used in [AVV6] to summarize the key
procedures needed for computer experiments, function tabulation,
and graphing. Although computation of most special functions is
 in principle ``well known,'' in practice 
finding algorithms requires much work and occasionally one has
to implement algorithms. Mathematica contains as built-ins
many of the functions we need.

A preliminary version of the software and the manual are available.
The software runs on both PC and Unix machines.

\bigskip
%\noindent{\bf 
\begin{nonsec} \label{3.63b}
Newton algorithm for $\mu^{-1}(y), y > {\pi \over 2}$ \end{nonsec}
Set $x_0 = 1/\mbox{cosh}\, y$ and
$$ x_{n+1}=x_n -{\mu(x_n)-y \over \mu'(x_n)} = 
x_n -{(\mu(x_n)-y)(x_n - x_n^3) \over AG(1, x'_n)^2} .$$

In practice, this algorithm always converges, but the proof of
convergence is missing.

\bigskip

%\noindent{\bf

\begin{nonsec} \label{3.64}  Open problems (for numerical analysis
students)
\end{nonsec} {

(a) Does the above Newton algorithm converge to $\mu^{-1}(y)$ for $y > \pi/2$?

(b) Is it true that $x_n < x_{n+1} < 1$ for all $n$ if $y > \pi$?
\noindent
%Note that (b)$ \Rightarrow$ (a). }
%$\K(r) = {\pi \over 2} F(1/2,1/2;1;r^2)$

\bigskip

D. Partyka \cite{Par1}-\cite{Par2} has also devised algorithms for computing
functions
related to $\mu^{-1}(y) .$

\bigskip

\begin{nonsec} \label{3.65}  Remark
\end{nonsec} { The open problems \ref{3.54} and \ref{3.64} have been
recently studied in \cite{AVV7}.}

\bigskip
{\bf Acknowledgements.} Many of my colleagues and friends have provided
helpful remarks on the earlier versions of this survey. I want to
thank them all. In particular, I am indebted to G. D. Anderson, J. Luukkainen,
M. K. Vamanamurthy, and P. Alestalo. I am also indebted to 
M. Kolehmainen, M. Nikunen, and P. Seittenranta, who made the pictures.

\bigskip

% Tulosta pienella fontilla

%POISTA KOMMENTTI SEUR: RIVIN ALUSTA JOS AJAT ITSENAISENA TIEDOSTONA
%\documentstyle[11pt,amssym,fullpage]{article}
%fullpage,twoside,amssym,samy
%\pagestyle{myheadings}

\def\R{\roman{Re}}
\def\I{\roman{Im}}
\def\mes{\mbox{\rm mes}}
\def\Rm{{{\Bbb R}^m}}
\def\Rn{{{\Bbb R}^n}}

%POISTA KOMMENTTI SEUR: RIVIN ALUSTA JOS AJAT ITSENAISENA TIEDOSTONA
%\begin{document}
%\addcontentsline{toc}{part}{\protect\numberline{}{REFERENCES}}

%\tiny
%\scriptsize

%\small

\noindent
%\today

\bigskip
\bigskip

%\centerline{\large{\bf \sc{Lublin Lectures by M. Vuorinen, August 1996}}}

\end{document}